\documentclass[aos,MSNbibl,nameyear,dvips]{arximspdf}
\usepackage{dcolumn}

\doi{10.1214/11-AOS877}
\volume{39}
\issue{2}
\pubyear{2011}
\firstpage{1310}
\lastpage{1333}

\makeatletter

  \let\sv@tabnotetext\tabnotetext
  \let\sv@tabnotemark@fmt\tabnotemark@fmt
   \long\def\legend#1{{\let\tabnote@indent\leavevmode\sv@tabnotetext[]{}{#1}}}

\newcolumntype{d}[1]{D{.}{.}{#1}}
\newtheorem{theo}{Theorem}
\newtheorem{lema}{Lemma}
\newtheorem{coro}{Corollary}
\newproclaim{example}{Example}
\newcommand{\eqref}[1]{(\ref{#1})}
\makeatother

\begin{document}
\begin{frontmatter}

\title{On construction of optimal mixed-level supersaturated designs}
\runtitle{Construction of mixed-level supersaturated designs}

\begin{aug}
\author[A]{\fnms{Fasheng} \snm{Sun}\thanksref{tnx}\ead[label=e1]{sfxsfx2001@gmail.com}},
\author[B]{\fnms{Dennis K. J.} \snm{Lin}\ead[label=e2]{DKL5@psu.edu}}
\and
\author[C]{\fnms{Min-Qian} \snm{Liu}\corref{}\thanksref{tnx}\ead[label=e3]{mqliu@nankai.edu.cn}}
\runauthor{F. Sun, D. K. J. Lin and M.-Q. Liu}
\affiliation{Northeast Normal University,   Pennsylvania State
University and~Nankai~University}
\address[A]{F. Sun\\
Department of Statistics\\
KLAS and School of Mathematics\\
\quad  and Statistics\\
Northeast Normal University\\
Changchun 130024\\
China\\
\printead{e1}} 
\address[B]{D. K. J. Lin\\
Department of Statistics\\
Pennsylvania State University\\
University Park,
Pennsylvania 16802\\
USA\\
\printead{e2}}
\address[C]{M.-Q. Liu\\
Department of Statistics\\
School of Mathematical Sciences\\
\quad  and LPMC\\
Nankai University\\
Tianjin 300071\\
China\\
\printead{e3}}
\end{aug}
\thankstext{tnx}{Supported by NNSF of China Grant 10971107 and Program
for New Century Excellent Talents in University (NCET-07-0454) of China.}

\received{\smonth{7} \syear{2010}}
\revised{\smonth{12} \syear{2010}}

%
\begin{abstract}
Supersaturated design (SSD) has received much recent interest
because of its potential in factor screening experiments.
In this paper, we provide equivalent conditions for two columns to
be fully aliased and consequently propose methods for
constructing $E(f_{\mathrm{NOD}})$- and $\chi^2$-optimal mixed-level SSDs
without fully aliased columns, via equidistant designs and
difference matrices. The methods can be easily performed and many
new optimal mixed-level SSDs have been obtained.
Furthermore, it is proved that the
nonorthogonality between columns of the resulting design is well
controlled by the source designs.
A rather complete list of newly generated optimal mixed-level SSDs are
tabulated for practical
use.
\end{abstract}

%

\end{frontmatter}

\section{Introduction}\label{ssec1}

The supersaturated design (SSD) is a factorial design in
which the number of runs is not sufficient to estimate all the main
effects. Such designs are useful when the experiment is expensive,
the number of factors is large, and only a few significant factors
need to be identified in a relatively small number of experimental
runs. \citet{BooCox62} first examined these designs
systematically and proposed the $E(s^{2})$ criterion.
However, such designs were not further studied until the appearance of the
work by Lin (\citeyear{Lin93,Lin95}), \citet{Wu93},
\citet{TanWu97} and \citet{CheTan01}. Research on
mixed-level SSDs includes the early work
by Fang, Lin and Liu (\citeyear{FanLinLiu,FanLinLiu03}) who
proposed the $E(f_{\mathrm{NOD}})$
criterion and the FSOA method for constructing mixed-level SSDs, and
work by \citet{YamMat02} and \citet{YamLin02} who used
$\chi^{2}$ to evaluate mixed-level SSDs. Recent work on
mixed-level SSDs includes \citet{Xu03}, \citet{Fanetal04N1},
\citet{LiLiuZha04}, \citet{XuWu05}, \citet{KouMan05},
\citet{LiuFanHic06},
Yamada et al. (\citeyear{Yametal06}),
\citet{AiFanHe07}, \citet{Tanetal07}, Chen and Liu (\citeyear
{CheLiu08N1,CheLiu08N2}),
\citet{LiuLin09}, \citet{LiuCai09} and \citet{LiuZha09}.

This paper proposes some methods for constructing
$E(f_{\mathrm{NOD}})$- and $\chi^2$-optimal mixed-level SSDs without fully
aliased columns, and with a control on the nonorthogonality.
A large number of optimal designs is obtained.
The remainder of this paper is organized as
follows. Section \ref{ssec2} provides relevant notation and definitions.
In Section \ref{ssec3}, we propose the general construction methods
for mixed-level SSDs along with illustrative examples. Discussions
on the nonorthogonality of the resulting designs are given in
Section \ref{ssec4}. In Section \ref{ssec5},
a review of the existing methods for mixed-level SSDs
and comparisons with the current methods are made,
and some concluding remarks are provided.
For coherence of presentation, all proofs are placed in Appendix~\ref{apdxa} and
newly constructed designs are tabulated in Appendix \hyperref[apdxb]{B}.

\section{Preliminaries}\label{ssec2}

A mixed-level design that has $n$ runs and $m$ factors
with $q_1, \ldots , q_m$ levels, respectively, is denoted by $F(n,
q_1\cdots q_m)$. When\break $\sum_{j=1}^m(q_j-1)=n-1$, the design is
called a saturated design, and when $\sum_{j=1}^m(q_j-1)>n-1$, the
design is called a supersaturated design (SSD).
An $F(n, q_1\cdots q_m)$ can be expressed as an $n \times m$ matrix
$F=(f_{ij})$. When some $q_j$'s are equal, we use the notation $F(n,
q_1^{r_1}\cdots q_l^{r_l})$ indicating $r_i$ factors having $q_i$ levels,
$i=1, \ldots, l$. If all the $q_j$'s are equal, the
design is said to be symmetrical and denoted by $F(n, q^m)$. Let
$f_i$ be the $i$th row of an $F(n, q_1\cdots q_m)$ and $f^j$ be the
$j$th column which takes values from a set of $q_j$ symbols
$\{0,\ldots,q_{j}-1\}$. If each column $f^j$ is balanced, that is, it
contains the $q_j$ symbols equally often, then we say $F$ is a
balanced design. Throughout this paper, we only consider balanced
designs. Two columns are called fully aliased if one column can be
obtained from the other by permuting levels;
and called orthogonal if all possible
level-combinations for these two columns appear equal number of times.
An $F(n, q_1\cdots q_m)$ is called an orthogonal array of strength two,
denoted by $L_n(q^m)$ for the symmetrical case, if all pairs of columns
of this design are orthogonal.

The set of residues modulo a prime number $p, \{0, 1, \ldots,
p-1\}$, forms a field of $p$ elements under addition and
multiplication modulo $p$, which is called a~Galois field and
denoted by $\operatorname{GF}(p)$. Note that the order of a Galois field must be a
prime power. A Galois field of order $q=p^{u}$ for any prime $p$
and any positive integer $u$ can be obtained as follows. Let
$g(x)=b_0+b_1x+\cdots+ b_ux^{u}$ be an irreducible polynomial of
degree $u$, where $b_j\in \operatorname{GF}(p)$ and $b_u = 1$. Then the set of all
polynomials of degree $u-1$ or lower,
$\{a_0+a_1x+\cdots+a_{u-1}x^{u-1}|a_j\in \operatorname{GF}(p)\},$ is a Galois field
$\operatorname{GF}(q)$ of order $q=p^{u}$ under addition and multiplication of
polynomials modulo $g(x)$. For any polynomial $f(x)$ with
coefficients from $\operatorname{GF}(p)$, there exist unique polynomials $q(x)$ and
$r(x)$ such that $f(x) = q(x)g(x)+r(x),$ where the degree of $r(x)$
is lower than $u$. This $r(x)$ is the residue of $f(x)$ modulo
$g(x)$, which is usually written as $f(x) = r(x) (\operatorname{mod} g(x)$).

A difference matrix, denoted by $D(rq, c, q)$, is an $rq\times c$
array with entries from a finite Abelian group $(\mathcal{A}, +)$
with $q$ elements such that each element of $\mathcal{A}$ appears
equally often in the vector of difference between any two columns of
the array [\citet{BosBus52}]. Note that if $A$ is an
$L_{rq}(q^c)$, then it is also a difference
matrix. A difference matrix $D(rq, c, q)$ with $c>1$ is said to be normalized,
denoted by $\operatorname{ND}(rq, c, q)$, if its first column consists of all
zeros.\vadjust{\goodbreak} In fact, for any difference matrix $D$, if we subtract the
first column from any column, then we can obtain a normalized
difference matrix.\looseness=1

For a scalar $a$ and a matrix $A$, let $a+A$ denote the
element-wise sum of $a$ and $A$. For any two matrices
$A=(a_{ij})$ of order $r\times s$ and $B$ of order $u\times v$,
their {Kronecker sum} and {Kronecker product} are defined to be\vspace{3pt}
\[
A\oplus B= \pmatrix{\displaystyle
a_{11}+B&\cdots&a_{1s}+B\cr\displaystyle
\cdots&\cdots&\cdots\cr\displaystyle
a_{r1}+B&\cdots&a_{rs}+B
}
   \quad \mbox{and} \quad  A\otimes B= \pmatrix{\displaystyle
a_{11}B&\cdots&a_{1s}B\cr\displaystyle
\cdots&\cdots&\cdots\cr\displaystyle
a_{r1}B&\cdots&a_{rs}B},\vspace{3pt}
\]
respectively.
Here, we use ``$+_{\mathcal{A}}$'' and ``$\oplus_{\mathcal{A}}$'' to denote
the sum and Kronecker sum defined on $\mathcal{A}$,
respectively.

For a design $F=(f_{ij})_{n\times m}$, let\vspace{3pt}
\begin{eqnarray*}
\lambda_{ij}(F)=\sum_{k=1}^{m} \delta_{ij}^{(k)}  \quad \mbox{and} \quad
\omega_{ij}(F)=\sum_{k=1}^{m} q_k\delta_{ij}^{(k)},
\end{eqnarray*}
where $\delta_{ij}^{(k)}=1$ if $f_{ik}=f_{jk}$, and 0 otherwise.
Then $\lambda_{ij}(F)$ and $\omega_{ij}(F)$ are called the
{coincidence number} and {natural weighted coincidence number}
between rows $f_i$ and $f_j$, respectively. A design with equal
coincidence numbers between different rows is called an {equidistant
design}. From \citet{MukWu95}, a saturated $L_n(q^m)$ is an
equidistant design with
%
\begin{equation}\label{sMW}
\lambda_{ij}(F)=\frac{m-1}{q}  \quad \mbox{and} \quad  \omega_{ij}(F)=m-1  \qquad \mbox{for } i\neq j.
\end{equation}

The $E(f_{\mathrm{NOD}})$ criterion proposed by Fang, Lin and Liu (\citeyear
{FanLinLiu,FanLinLiu03}) is defined to minimize
\[
E(f_{\mathrm{NOD}})= \frac{2}{m(m-1)}\sum_{1\leq i<j\leq m}f_{\mathrm{NOD}}(f^i,
f^j),
\]
 where
 \[
f_{\mathrm{NOD}}(f^i, f^j)= \sum^{q_i-1}_{a=0} \sum^{q_j-1}_{b=0} \biggl(n_{ab}(f^i,
f^j)-\frac{n}{q_iq_j} \biggr)^{2},
\]
$n_{ab}(f^i, f^j)$ is the number of $(a, b)$-pairs in $(f^i, f^j)$,
and $n/(q_i q_j)$ stands for the average frequency of level-combinations
in $(f^i, f^j)$. Here, the subscript ``NOD'' stands for \textit{nonorthogonality of the design}.
The $f_{\mathrm{NOD}}(f^i, f^j)$ value gives a nonorthogonality measure for
$(f^i, f^j)$, and
columns $f^i$ and $f^j$ are orthogonal if and only if $f_{\mathrm{NOD}}(f^i, f^j)=0$.
It is obvious that $F$ is an orthogonal array if and only if
$E(f_{\mathrm{NOD}})=0$, that is,
$f_{\mathrm{NOD}}(f^i, f^j)=0$ for all $i, j=1,\ldots, m, i\neq j$. Thus
$E(f_{\mathrm{NOD}})$
measures the average nonorthogonality among the columns of $F$.

Another criterion that is to be minimized was defined by\vspace*{1pt} \citet{YamLin99} and
\citet{YamMat02} as
$\chi^2(F)=\sum_{1\leq i<j\leq m}q_iq_jf_{\mathrm{NOD}}(f^i\!,\break f^j)/n$.
Obviously,
$E(f_{\mathrm{NOD}})$ and $\chi^2(F)$ are equivalent in
the symmetrical case. Here, we adopt both $E(f_{\mathrm{NOD}})$ and
$\chi^2(F)$ to evaluate the newly constructed SSDs. There are also
some other criteria for assessing mixed-level SSDs [see, e.g., \citet{LiuLin09} for a general review].

The following results, regarding the $E(f_{\mathrm{NOD}})$ and $\chi^2(F)$
optimality criteria of a design, will be needed for our construction
methods.\vspace{3pt}

\begin{lema}\label{slem1}
\textup{(a) [\citet{Fanetal04N1}].} If the difference among
all coincidence numbers between different rows of design $F$ does
not exceed one, then $F$ is $E(f_{\mathrm{NOD}})$-optimal.\vspace{3pt}

\textup{(b) [\citet{LiLiuZha04}; \citet{LiuFanHic06}].}
If the natural weighted coincidence numbers between different rows
of design $F$ take at most two nearest values, then $F$ is $\chi^2$-optimal.
\end{lema}

\section{Proposed construction methods}\label{ssec3}

In this section, we first provide some equivalent conditions for two
columns to
be fully aliased, then propose methods for constructing
$E(f_{\mathrm{NOD}})$- and $\chi^2$-optimal SSDs, and finally study the
properties of the resulting designs.

\subsection{Equivalent conditions for two columns to
be fully aliased}

An\break $E(f_{\mathrm{NOD}})$- or $\chi^2$-optimal SSD may contain fully aliased
columns, which is
undesirable. Let matrix $X_j=(x^j_{st})$ of order $n\times q_j$ be the
induced matrix
[\citet{Fanetal04N1}] of the $j$th column of an $F(n, q_1\cdots
q_m)$, that is, $x^j_{st}=1$ if the $s$th element in the $j$th column
is $t-1$, otherwise 0, for $s=1, \ldots, n, t=1, \ldots, q_j$ and
$j=1, \ldots, m.$ The following theorem presents theoretical
results concerning the column aliasing that will be used in the
construction methods.\vspace{3pt}

\begin{theo}\label{slem0}
Suppose $X_j=(x^j_{st})$ is the induced matrix of a balanced \break column
$f^j=(f_{1j}, \ldots, f_{n_jj})'$ with $q_j$ levels, $j=1, \ldots,
4$, and $n_1=n_3,\break n_2=n_4$.
\begin{enumerate}[(a)]
\item[(a)] For $q_1=q_2=q_3=q_4=q$ and ${\mathcal{A}}=\{0,
\ldots,
q-1\}$:
\begin{longlist}[(iii)]
\item[(i)] $f^1$ and $f^3$ are fully aliased if and only if
$X_1X'_1=X_3X'_3$;\vspace*{2pt}
\item[(ii)] the induced matrix of $f^1\oplus_{\mathcal{A}}f^2$ is\vspace*{2pt}
$[(X_2P_{f_{11}})', \ldots, (X_2P_{f_{n_11}})']'
=(X_1\otimes X_2)P$, where $P= (P'_0, \ldots,
P'_{q-1} )'$ and $P_{i}$ is a permutation matrix
defined by
\[
i+_{{\mathcal{A}}}(0, \ldots, q-1)=(0, \ldots, q-1)P'_{i},\qquad i=0,
\ldots, q-1;\vspace*{3pt}
\]
\item[(iii)] if $f^1\oplus_{\mathcal{A}}f^2$ and $f^3\oplus
_{\mathcal{A}}f^4$ are fully aliased,
then $f^1$ is fully aliased with $f^3$ and $f^2$ is fully aliased
with $f^4$.
\end{longlist}
\item[(b)]
\begin{longlist}[(iii)]
\item[(i)] The induced matrix of the $q_1q_2$-level column\vspace*{1pt}
$q_2(f^1-\frac{q_1-1}{2})\oplus(
f^2-\frac{q_2-1}{2})+\frac{q_1q_2-1}{2}$ is $X_1\otimes X_2$;\vspace*{1pt}
\item[(ii)] columns $q_2(f^1-\frac{q_1-1}{2})\oplus(f^2-\frac
{q_2-1}{2})+\frac{q_1q_2-1}{2}$
and\vspace*{2pt} $q_4(f^3-\frac{q_3-1}{2})\oplus(f^4-\frac{q_4-1}{2})+\frac
{q_3q_4-1}{2}$ are fully
aliased if and only\vspace*{3pt} if $f^1$ is fully aliased with $f^3$ and $f^2$ is fully
aliased with $f^4$;\vspace*{2pt}
\item[(iii)] for\vspace*{2pt} $q_3=q_4=q$, $q_2(f^1-\frac{q_1-1}{2})\oplus
(f^2-\frac{q_2-1}{2})+\frac{q_1q_2-1}{2}$
and $f^3\oplus_{\mathcal{A}}f^4$ are not fully aliased in any case.\vspace*{2pt}
\end{longlist}
\end{enumerate}
\end{theo}

\subsection{Construction of optimal symmetrical SSDs}

We next present the~me\-thods for constructing $E(f_{\mathrm{NOD}})$-
and $\chi^2$-optimal SSDs without fully aliased columns.

\begin{theo}\label{sthm1}
Let $D$ be an $\operatorname{ND}(rq, c, q)$ defined on an Abelian group
$\mathcal{A}=\{0, \ldots, q-1\}$ without identical rows, $F$ be an
$F(n, q^{m})$ without fully aliased columns and with constant
coincidence numbers, say $\lambda$, between its different rows,
then:
\begin{longlist}[(b)]
\item[(a)] $F\oplus_{\mathcal{A}}D'$ is an $F(cn, q^{rqm})$ with
two different values of coincidence numbers, $mr$ and $\lambda
rq$;
\item[(b)] $F\oplus_{\mathcal{A}}D'$ has no fully aliased columns.
\end{longlist}
\end{theo}

From Lemma \ref{slem1}, if $|mr-\lambda rq|\leq1$, then
$F\oplus_{\mathcal{A}} D'$ is both $E(f_{\mathrm{NOD}})$- and
$\chi^2$-optimal. The following corollary can be directly obtained from
Lemma~\ref{slem1}, Theorem~\ref{sthm1}, and equation \eqref{sMW}.

\begin{coro}\label{scoro1}
Let $F$ be a saturated $L_n(q^m)$ and $D$ be an $\operatorname{ND}(q, c, q)$
without identical rows. Then $F\oplus_{\mathcal{A}} D'$ is an
$F(cn, q^{mq})$ without fully aliased columns and with two different
values of coincidence numbers, $m$ and $m-1$, and thus is both
$E(f_{\mathrm{NOD}})$- and $\chi^2$-optimal.
\end{coro}

From Hedayat, Slone and Stufken (\citeyear{HedSloStu99}), there exist an
$L_{n}(q^{m})$ with $n=q^t$ and $m=(n-1)/(q-1)$ and an $\operatorname{ND}(q,
q, q)$ without identical rows for any prime power $q$, thus optimal $F(cq^t,
q^{(q^{t+1}-q)/(q-1)})$\vspace*{2pt} designs with coincidence numbers
$(q^t-1)/(q-1)-1$ or
$(q^t-1)/(q-1)$ can be constructed from Corollary \ref{scoro1},
where $c$ is a positive integer and $c<q$.

\begin{example}\label{se1}
Let $F$ be an $L_9(3^4)$ and $D$ be an $\operatorname{ND}(3, 2, 3)$ (cf. Table
\ref{stb1}), then $F\oplus_{\mathcal{A}} D'$ is an $F(18, 3^{12})$
with coincidence numbers $4$ and $3$ as listed in Table \ref{stb3},
where $\mathcal{A}=\operatorname{GF}(3)$. This new design is an $E(f_{\mathrm{NOD}})$-
and $\chi^2$-optimal SSD without fully aliased columns.
\end{example}

\begin{table}
\vspace{4pt}
\tabcolsep=0pt
\caption{$F$ and $D$ in Example \protect\ref{se1}}\label{stb1}
\begin{tabular*}{170pt}{@{\extracolsep{\fill}}lcccccccc@{\hspace*{10pt}}cc@{}}
\hline
\multicolumn{9}{@{}c@{\hspace*{10pt}}}{$\bolds{F'}$}&\multicolumn{2}{c@{}}{$\bolds D$}\\
\hline
0&0&0&1&1&1&2&2&2 & 0&0 \\
0&1&2&0&1&2&0&1&2 & 0&1 \\
0&1&2&1&2&0&2&0&1 & 0&2 \\
0&2&1&1&0&2&2&1&0 & & \\
\hline
\end{tabular*}
\end{table}
%
%

\begin{table}[b]
\caption{The $F(18, 3^{12})$ constructed in Example \protect\ref{se1}}\label{stb3}
\begin{tabular}{@{}lccccccccccc@{}}
 \hline
\multicolumn{12}{@{}c@{}}{$\bolds{F\oplus_{\mathcal{A}} D'}$}\\
\hline
0&0&0&0&0&0&0&0&0&0&0&0\\
0&1&2&0&1&2&0&1&2&0&1&2\\
0&0&0&1&1&1&1&1&1&2&2&2\\
0&1&2&1&2&0&1&2&0&2&0&1\\
0&0&0&2&2&2&2&2&2&1&1&1\\
0&1&2&2&0&1&2&0&1&1&2&0\\
1&1&1&0&0&0&1&1&1&1&1&1\\
1&2&0&0&1&2&1&2&0&1&2&0\\
1&1&1&1&1&1&2&2&2&0&0&0\\
1&2&0&1&2&0&2&0&1&0&1&2\\
1&1&1&2&2&2&0&0&0&2&2&2\\
1&2&0&2&0&1&0&1&2&2&0&1\\
2&2&2&0&0&0&2&2&2&2&2&2\\
2&0&1&0&1&2&2&0&1&2&0&1\\
2&2&2&1&1&1&0&0&0&1&1&1\\
2&0&1&1&2&0&0&1&2&1&2&0\\
2&2&2&2&2&2&1&1&1&0&0&0\\
2&0&1&2&0&1&1&2&0&0&1&2\\
\hline
\end{tabular}
\end{table}
%
%
\subsection{Construction of optimal SSDs with two different level sizes}

Based on Lemma \ref{slem1} and Theorem \ref{sthm1},
the following theorem can be obtained.
\begin{theo}\label{sthm2}
Let $F_i$ be an $F(n_i, q^{m_i}_i)$ with constant coincidence
numbers $\lambda_i$, and no full aliased columns, $i=1,2$. Let
$D$ be an $\operatorname{ND}(rq_1, n_2, q_1)$ defined on Abelian group
$\mathcal{A}_1=\{0, \ldots, q_1-1\}$ without\vadjust{\goodbreak} identical rows. Then
$F=(F_1\oplus_{\mathcal{A}_1}D',  0_{n_1}\oplus F_2)$ is an
$F(n_1n_2, q^{rm_1q_1}_1q^{m_2}_2)$ without full aliased columns.
Furthermore:
\begin{longlist}[(b)]
\item[(a)] if $|(\lambda_2+rm_1)-(m_2+\lambda_1rq_1)|\leq1$,
then $F$
is $E(f_{\mathrm{NOD}})$-optimal;
\item[(b)] if $q_2\lambda_2+q_1rm_1=q_2m_2+\lambda_1rq^2_1$, then $F$
is $\chi^2$-optimal.
\end{longlist}
\end{theo}

Next, let us consider two illustrative examples for Theorem \ref{sthm2}.

\begin{example}\label{se2}
Let $F_1$ be an $L_4(2^3)$, $F_2$ be the $E(f_{\mathrm{NOD}})$-optimal $F(6,
3^{5})$ obtained by \citet{FanGeLiu04} and $D$ be an $\operatorname{ND}(8, 6,
2)$ without identical rows
obtained from an $L_8(2^7)$ based on $\mathcal{A}=\operatorname{GF}(2)$. Then
$\lambda_1=\lambda_2=1, q_1=2, q_2=3, r=4, m_1=3 \mbox{ and } m_2=5$
which satisfy the condition that
$\lambda_2+rm_1=m_2+\lambda_1rq_1=13$, thus $(F_1\oplus_{\mathcal
{A}}D', 0_{4}\oplus F_2)$ is an $E(f_{\mathrm{NOD}})$-optimal\vadjust{\goodbreak} $F(24,
2^{24}3^{5})$ with constant coincidence numbers $13$. The source
designs and resulting design are listed in Tables \ref{stb4} and
\ref{stb5}, respectively.
%
\begin{table}
\caption{$F_1, F_2$ and $D$ in Example \protect\ref{se2}}\label{stb4}
\tabcolsep=0pt
\begin{tabular*}{253pt}{@{\extracolsep{\fill}}lcc@{\hspace*{10pt}}ccccc@{\hspace*{10pt}}cccccccc@{}} \hline
\multicolumn{3}{@{}c@{\hspace*{10pt}}}{$\bolds{F_1}$}&\multicolumn{5}{c@{\hspace*{10pt}}}{$\bolds{F_2}$}&\multicolumn
{8}{c@{}}{$\bolds{D'}$}\\
\hline
0 & 0 & 0 & 0& 0 &0& 0 &0 & 0 & 0 & 0 & 0 & 0& 0 & 0& 0 \\
0 & 1 & 1 & 0& 1 &1& 1 &1 & 0 & 0 & 1 & 1 & 1& 1 & 0& 0 \\
1 & 0 & 1 & 1& 0 &2& 2 &1 & 0 & 1 & 0 & 1 & 1& 0 & 1& 0 \\
1 & 1 & 0 & 1& 2 &0& 1 &2 & 0 & 1 & 0 & 1 & 0& 1 & 0& 1 \\
& & & 2& 1 &2& 0 &2 & 0 & 0 & 1 & 1 & 0& 0 & 1& 1 \\
& & & 2& 2 &1& 2 &0 & 0 & 1 & 1 & 0 & 0& 1 & 1& 0 \\
\hline
\end{tabular*}
\end{table}
%

%
\begin{table}
\tabcolsep=0pt
\caption{The $F(24, 2^{24}3^{5})$ constructed in Example
\protect\ref{se2}} \label{stb5}
\begin{tabular*}{\textwidth}{@{\extracolsep{4in minus 4in}}lccccccccccccccccccccccc@{\hspace*{10pt}}ccccc@{}}
\hline
\multicolumn{24}{@{}c@{\hspace*{10pt}}}{$\bolds{F_1\oplus_{\mathcal{A}}D'}$}&\multicolumn
{5}{c@{}}{$\bolds{0_{4}\oplus F_2}$}\\
\hline
0 & 0 & 0 & 0 & 0& 0& 0& 0& 0 & 0& 0& 0& 0& 0 & 0& 0 & 0 & 0 & 0& 0 & 0
& 0 & 0 & 0 & 0 & 0& 0& 0 & 0 \\
0 & 0 & 1 & 1 & 1& 1& 0& 0& 0 & 0& 1& 1& 1& 1 & 0& 0 & 0 & 0 & 1& 1 & 1
& 1 & 0 & 0 & 0 & 1& 1& 1 & 1 \\
0 & 1 & 0 & 1 & 1& 0& 1& 0& 0 & 1& 0& 1& 1& 0 & 1& 0 & 0 & 1 & 0& 1 & 1
& 0 & 1 & 0 & 1 & 0& 2& 2 & 1 \\
0 & 1 & 0 & 1 & 0& 1& 0& 1& 0 & 1& 0& 1& 0& 1 & 0& 1 & 0 & 1 & 0& 1 & 0
& 1 & 0 & 1 & 1 & 2& 0& 1 & 2 \\
0 & 0 & 1 & 1 & 0& 0& 1& 1& 0 & 0& 1& 1& 0& 0 & 1& 1 & 0 & 0 & 1& 1 & 0
& 0 & 1 & 1 & 2 & 1& 2& 0 & 2 \\
0 & 1 & 1 & 0 & 0& 1& 1& 0& 0 & 1& 1& 0& 0& 1 & 1& 0 & 0 & 1 & 1& 0 & 0
& 1 & 1 & 0 & 2 & 2& 1& 2 & 0 \\
0 & 0 & 0 & 0 & 0& 0& 0& 0& 1 & 1& 1& 1& 1& 1 & 1& 1 & 1 & 1 & 1& 1 & 1
& 1 & 1 & 1 & 0 & 0& 0& 0 & 0 \\
0 & 0 & 1 & 1 & 1& 1& 0& 0& 1 & 1& 0& 0& 0& 0 & 1& 1 & 1 & 1 & 0& 0 & 0
& 0 & 1 & 1 & 0 & 1& 1& 1 & 1 \\
0 & 1 & 0 & 1 & 1& 0& 1& 0& 1 & 0& 1& 0& 0& 1 & 0& 1 & 1 & 0 & 1& 0 & 0
& 1 & 0 & 1 & 1 & 0& 2& 2 & 1 \\
0 & 1 & 0 & 1 & 0& 1& 0& 1& 1 & 0& 1& 0& 1& 0 & 1& 0 & 1 & 0 & 1& 0 & 1
& 0 & 1 & 0 & 1 & 2& 0& 1 & 2 \\
0 & 0 & 1 & 1 & 0& 0& 1& 1& 1 & 1& 0& 0& 1& 1 & 0& 0 & 1 & 1 & 0& 0 & 1
& 1 & 0 & 0 & 2 & 1& 2& 0 & 2 \\
0 & 1 & 1 & 0 & 0& 1& 1& 0& 1 & 0& 0& 1& 1& 0 & 0& 1 & 1 & 0 & 0& 1 & 1
& 0 & 0 & 1 & 2 & 2& 1& 2 & 0 \\
1 & 1 & 1 & 1 & 1& 1& 1& 1& 0 & 0& 0& 0& 0& 0 & 0& 0 & 1 & 1 & 1& 1 & 1
& 1 & 1 & 1 & 0 & 0& 0& 0 & 0 \\
1 & 1 & 0 & 0 & 0& 0& 1& 1& 0 & 0& 1& 1& 1& 1 & 0& 0 & 1 & 1 & 0& 0 & 0
& 0 & 1 & 1 & 0 & 1& 1& 1 & 1 \\
1 & 0 & 1 & 0 & 0& 1& 0& 1& 0 & 1& 0& 1& 1& 0 & 1& 0 & 1 & 0 & 1& 0 & 0
& 1 & 0 & 1 & 1 & 0& 2& 2 & 1 \\
1 & 0 & 1 & 0 & 1& 0& 1& 0& 0 & 1& 0& 1& 0& 1 & 0& 1 & 1 & 0 & 1& 0 & 1
& 0 & 1 & 0 & 1 & 2& 0& 1 & 2 \\
1 & 1 & 0 & 0 & 1& 1& 0& 0& 0 & 0& 1& 1& 0& 0 & 1& 1 & 1 & 1 & 0& 0 & 1
& 1 & 0 & 0 & 2 & 1& 2& 0 & 2 \\
1 & 0 & 0 & 1 & 1& 0& 0& 1& 0 & 1& 1& 0& 0& 1 & 1& 0 & 1 & 0 & 0& 1 & 1
& 0 & 0 & 1 & 2 & 2& 1& 2 & 0 \\
1 & 1 & 1 & 1 & 1& 1& 1& 1& 1 & 1& 1& 1& 1& 1 & 1& 1 & 0 & 0 & 0& 0 & 0
& 0 & 0 & 0 & 0 & 0& 0& 0 & 0 \\
1 & 1 & 0 & 0 & 0& 0& 1& 1& 1 & 1& 0& 0& 0& 0 & 1& 1 & 0 & 0 & 1& 1 & 1
& 1 & 0 & 0 & 0 & 1& 1& 1 & 1 \\
1 & 0 & 1 & 0 & 0& 1& 0& 1& 1 & 0& 1& 0& 0& 1 & 0& 1 & 0 & 1 & 0& 1 & 1
& 0 & 1 & 0 & 1 & 0& 2& 2 & 1 \\
1 & 0 & 1 & 0 & 1& 0& 1& 0& 1 & 0& 1& 0& 1& 0 & 1& 0 & 0 & 1 & 0& 1 & 0
& 1 & 0 & 1 & 1 & 2& 0& 1 & 2 \\
1 & 1 & 0 & 0 & 1& 1& 0& 0& 1 & 1& 0& 0& 1& 1 & 0& 0 & 0 & 0 & 1& 1 & 0
& 0 & 1 & 1 & 2 & 1& 2& 0 & 2 \\
1 & 0 & 0 & 1 & 1& 0& 0& 1& 1 & 0& 0& 1& 1& 0 & 0& 1 & 0 & 1 & 1& 0 & 0
& 1 & 1 & 0 & 2 & 2& 1& 2 & 0 \\
\hline
\end{tabular*}
\end{table}
\end{example}

\begin{example}\label{se22}
Let $F_1$ be an $L_4(2^3)$, $F_2$ be the $F(6,
3^{10})$ obtained by \citet{GeoKou06} and $D$ be an $\operatorname{ND}(24, 6,
2)$\vadjust{\goodbreak} without
identical rows obtained from an $L_{24}(2^{23})$ based on $\mathcal
{A}=\operatorname{GF}(2)$. Then $\lambda_1=1,\break\lambda_2=2, q_1=2, q_2=3, r=12,
m_1=3 \mbox{ and } m_2=10$ which satisfy the condition that
$q_2\lambda_2+q_1rm_1=q_2m_2+\lambda_1rq^2_1=78$, thus
$(F_1\oplus_{\mathcal{A}}D', 0_{4}\oplus F_2)$ is a
$\chi^2(D)$-optimal $F(24, 2^{72}3^{10})$ with constant natural
weighted coincidence numbers~$78$.
Exact details are omitted here but available upon request.\looseness=1
\end{example}

\subsection{Construction of optimal SSDs with three different level sizes}

The next lemma is useful in the upcoming proposed construction method.
\begin{lema}\label{slem3}
Let $V=\{-\frac{q-1}{2}, -\frac{q-3}{2}, \ldots, \frac{q-3}{2},
\frac{q-1}{2}\}=\{0, \ldots, q-1\}-\frac{q-1}{2}$ and
$V_i=(i-\frac{p-1}{2})q+V, i=0, \ldots, p-1$, then $V_i\cap
V_j=\Phi$ for $i\neq j$ and $\bigcup^{p-1}_{i=0}V_i= \{-\frac{pq-1}{2}, -\frac{pq-3}{2}, \ldots, \frac{pq-3}{2},
\frac{pq-1}{2}\}= \{0, \ldots, pq-1\}-\frac{pq-1}{2}$, where
$\Phi$ is an empty set.
\end{lema}

From this lemma, we can obtain the following theorem in a
straightforward manner.
\begin{theo}\label{sth3}
Let $F_i$ be an $F(n_i, q^{m_i}_i)$ with constant coincidence
numbers $\lambda_i,$ $i=1, 2$, then
$q_2(F_1-\frac{q_1-1}{2})\oplus(F_2-\frac{q_2-1}{2})+\frac{q_1q_2-1}{2}$
is an $F(n_1n_2, (q_1q_2)^{m_1m_2})$ with three different values of
coincidence numbers $\lambda_1m_2, \lambda_2m_1$ and $
\lambda_1\lambda_2$.
\end{theo}

This theorem, along with Lemma \ref{slem0} and Theorem \ref{sthm1},
leads to the following theorem, which provides another
construction method of $E(f_{\mathrm{NOD}})$- and $\chi^2$-optimal SSDs.

\begin{theo}\label{sthm4}
Suppose $F_i$ is an $F(n_i, q^{m_i}_i)$ with constant coincidence
numbers $\lambda_i$ and no fully aliased columns, $i=1, \ldots, 4 $,
$D_3$ is an $\operatorname{ND}(r_3q_{3}, n_{2}, q_3)$ defined on Abelian group
$\mathcal{A}_3=\{0, \ldots, q_3-1\}$ without identical rows, $D_4$
is an $\operatorname{ND}(r_4q_{4}, n_{1}, q_4)$ defined on $\mathcal{A}_4=\{0,
\ldots, q_4-1\}$ without identical rows, and they satisfy \textup{(i)}
$n_1=n_3, n_2=n_4$; \textup{(ii)} the first rows of $F_3$ and $F_4$
consist of all zeros; \textup{(iii)}~there are no fully aliased columns
between $F_3$ and $D'_4$ or between $F_4$ and $D'_3$. Then
%
\begin{eqnarray}\label{s5}
F&=& \biggl[q_2 \biggl(F_1-\frac{q_1-1}{2} \biggr)\oplus \biggl(F_2-\frac
{q_2-1}{2} \biggr)+\frac{q_1q_2-1}{2},\nonumber
\\[-8pt]
\\[-8pt]
&&\hspace*{120pt}F_3\oplus_{\mathcal{A}_{3}}D'_3, D'_4\oplus_{\mathcal{A}_{4}}F_4 \biggr]
\nonumber
\end{eqnarray}
is an $F(n_1n_2, (q_1q_2)^{m_1m_2}q^{m_3r_3q_3}_3q^{m_4r_4q_4}_4)$
without fully aliased columns
and:
\begin{longlist}[(b)]
\item[(a)] if the difference among three values
$\lambda_2m_1+r_3m_3+\lambda_4r_4q_4, \lambda_1\lambda
_2+r_3m_3+r_4m_4$ and
$\lambda_1m_2+\lambda_3r_3q_3+r_4m_4$ does not exceed one, then $F$
is $E(f_{\mathrm{NOD}})$-optimal;
\item[(b)] if $q_1q_2\lambda_2m_1+q_3r_3m_3+\lambda_4r_4q^2_4=
q_1q_2\lambda_1\lambda_2+q_3r_3m_3+q_4r_4m_4=\break q_1q_2\lambda_1m_2+\lambda
_3r_3q^2_3+q_4r_4m_4$, then $F$ is $\chi^2$-optimal.
\end{longlist}
\end{theo}

The following two examples serve as illustrations of the
construction method in Theorem \ref{sthm4}.

\begin{example}
Let $F_1$ and $F_3$ be two $L_4(2^3)$'s; $F_2$ be the $F(6, 2^{10})$
obtained by \citet{LiuZha00};
$F_4$ be the $F(6, 3^5)$ obtained by \citet{FanGeLiu04};
$D_3$ be an $\operatorname{ND}(12, 6, 2)$ without identical rows obtained from an
$L_{12}(2^{11})$; $D_4$ be an $\operatorname{ND}(12, 4, 3)$ without identical rows;
$\mathcal{A}_{3}=\operatorname{GF}(2)$ and $\mathcal{A}_{4}=\operatorname{GF}(3)$. Suppose the
first rows of $F_3$ and $F_4$ consist of all zeros. Then based on
Theorem \ref{sthm4}, $\lambda_1=\lambda_3=\lambda_4=1, \lambda_2=4,
m_1=m_3=3, m_2=10, m_4=5, q_1=q_2=q_3=2, q_4=3, r_3=6, r_4=4$ and\vspace*{1pt}
$\lambda_2m_1+r_3m_3+\lambda_4r_4q_4=\lambda_1\lambda
_2+r_3m_3+r_4m_4=\lambda_1m_2+\lambda_3r_3q_3+r_4m_4=42$.
Thus, from \eqref{s5}, we obtain an $E(f_{\mathrm{NOD}})$-optimal $F(24,
4^{30}2^{36}3^{60})$ with constant coincidence numbers $42$ and no
fully aliased columns.
\end{example}

\begin{example}
Let $F_1$ and $F_3$ be two $L_4(2^3)$'s; both $F_2$ and $F_4$ be
the $F(6, 3^{5})$ obtained by \citet{FanGeLiu04}; $D_3$ be an
$\operatorname{ND}(24, 6, 2)$ without identical rows obtained from an
$L_{24}(2^{23})$ based on $\mathcal{A}_3=\operatorname{GF}(2)$ and $D_4$ be an
$\operatorname{ND}(6, 4, 3)$ without identical rows based on $\mathcal
{A}_4=\operatorname{GF}(3)$. Suppose the first rows of $F_3$ and $F_4$ consist of
all zeros. Then $\lambda_1=\lambda_2=\lambda_3=\lambda_4=1,
q_1=q_3=2, q_2=q_4=3, r_3=12, r_4=2, m_1=m_3=3, m_2=m_4=5$, which
satisfy the condition that
$q_1q_2\lambda_2m_1+q_3r_3m_3+\lambda_4r_4q^2_4=
q_1q_2\lambda_1\lambda_2+q_3r_3m_3+q_4r_4m_4=q_1q_2\lambda_1m_2+\lambda
_3r_3q^2_3+q_4r_4m_4=108$. Thus, the design constructed through \eqref
{s5} is a $\chi^2$-optimal $F(24, 6^{15}2^{72}3^{30})$ with constant
natural weighted
coincidence numbers $108$ and no fully aliased columns.
\end{example}

\section{Nonorthogonality of the resulting designs}\label{ssec4}

In the previous section, construction methods for $E(f_{\mathrm{NOD}})$- as
well as $\chi^2$-optimal SSDs without fully aliased columns are
provided. Full aliasing can be viewed as the extreme case of nonorthogonality.
In this section, we will investigate nonorthogonality, measured by
$f_{\mathrm{NOD}}$, of the resulting
designs, and show how it is controlled by the source
designs.

\begin{theo}\label{sthm6}
Suppose $f^i=(f_{1i}, \ldots, f_{n_ii})'$ is a $q_i$-level
balanced column with induced matrix $X_i$, $\mathcal{A}_i=\{0,
\ldots, q_i-1\}$, $i=1, \ldots, 4$, $n_1=n_3, n_2=n_4$. Let
$h_1=q_2(f^1-\frac{q_1-1}{2})\oplus(f^2-\frac{q_2-1}{2})+\frac{q_1q_2-1}{2}$
and
$h_2=q_4(f^3-\frac{q_3-1}{2})\oplus(f^4-\frac{q_4-1}{2})+\frac{q_3q_4-1}{2}$.
Then:
\begin{longlist}[(b)]
\item[(a)]
$f_{\mathrm{NOD}}(h_1, h_2)\vspace*{-2pt} =
f_{\mathrm{NOD}}(f^1,f^3)f_{\mathrm{NOD}}(f^2,f^4) +\frac{n^2_2}{q_2q_4}f_{\mathrm{NOD}}(f^1,f^3)+
\break\frac{n^2_1}{q_1q_3}f_{\mathrm{NOD}}(f^2, f^4);$\vspace*{1pt}

\item[(b)] if $q_1=q_2$, $q_3=q_4$, then
\begin{eqnarray*}
f_{\mathrm{NOD}}(f^1\oplus_{\mathcal{A}_1}f^2, f^3\oplus_{\mathcal{A}_3}f^4)
&\leq& q_1q_3f_{\mathrm{NOD}}(f^1, f^3)f_{\mathrm{NOD}}(f^2, f^4)\\
&&{}+\min \{n^2_2f_{\mathrm{NOD}}(f^1, f^3),n^2_1f_{\mathrm{NOD}}(f^2, f^4) \},
\end{eqnarray*}
where the equality holds if and only if $f^1$ is orthogonal to $f^3$ or
$f^2$ is orthogonal to $f^4$;

\item[(c)] if $q_1=q_2$, then
\begin{eqnarray*}
f_{\mathrm{NOD}}(f^1\oplus_{\mathcal{A}_1}f^2,
h_2) &\leq& q_1f_{\mathrm{NOD}}(f^1, f^3)f_{\mathrm{NOD}}(f^2, f^4)\\
&&{}+\min \{n^2_2/q_4f_{\mathrm{NOD}}(f^1, f^3), n^2_1/q_3f_{\mathrm{NOD}}(f^2,
f^4) \},
\end{eqnarray*}
where the equality holds if and only if $f^1$ is orthogonal to $f^3$ or
$f^2$ is orthogonal to $f^4$.
\end{longlist}
\end{theo}

Theorem \ref{sthm6} shows that the nonorthogonality measured by
$f_{\mathrm{NOD}}$ of the resulting designs is well controlled by the source
designs. If the source designs have small values of
$f_{\mathrm{NOD}}$, then the resulting design will also have small values
of $f_{\mathrm{NOD}}$. In particular, we have the following.
\begin{coro}
Suppose $f^i=(f_{1i}, \ldots, f_{n_ii})'$ is a $q_i$-level
balanced column with induced matrix $X_i$, $\mathcal{A}_i=\{0,
\ldots, q_i-1\}$, $i=1, \ldots, 4$, $n_1=n_3, n_2=n_4$. Then:
\begin{longlist}[(b)]
\item[(a)] if $f^1$ is orthogonal to $f^3$
and $f^2$ is orthogonal to $f^4$, then $q_2(f^1-\frac{q_1-1}{2})\oplus
(f^2-\frac{q_2-1}{2})+\frac{q_1q_2-1}{2}$
is orthogonal to
$q_4(f^3-\frac{q_3-1}{2})\oplus(f^4-\frac{q_4-1}{2})+\frac{q_3q_4-1}{2}$;

\item[(b)] if $q_1=q_2$, $q_3=q_4$ and $f^1$ is orthogonal to $f^3$
or $f^2$ is orthogonal to~$f^4$, then
$f^1\oplus_{\mathcal{A}_1}f^2$ is
orthogonal to $f^3\oplus_{\mathcal{A}_3}f^4$;

\item[(c)] if $q_1=q_2$ and $f^1$ is orthogonal to $f^3$
or $f^2$ is orthogonal to $f^4$, then $f^1\oplus_{\mathcal{A}_1}f^2$
is orthogonal to $q_4(f^3-\frac{q_3-1}{2})\oplus(f^4-\frac
{q_4-1}{2})+\frac{q_3q_4-1}{2}$.
\end{longlist}
\end{coro}

This corollary indicates that the orthogonality between columns of
the source design is maintained in the generated designs.

\section{Discussion and concluding remarks}\label{ssec5}

In this paper, we have presented some construction methods for $E(f_{\mathrm{NOD}})$-
and $\chi^2$-optimal SSDs. A review of the existing methods for
mixed-level SSDs
and comparisons with the current methods are summarized below.

(a) \citet{YamMat02} and \citet{YamLin02} proposed two
methods for
constructing mixed-level SSDs consisting of only two- and three-level columns
through computer searches. However, their
resulting designs have no theoretical support and typically are unable
to achieve
the lower bound of $\chi^2$-value.

(b) Fang, Lin and Liu (\citeyear{FanLinLiu,FanLinLiu03})
proposed an FSOA method for constructing
$E(f_{\mathrm{NOD}})$-optimal mixed-level SSDs from saturated orthogonal
arrays. \citet{LiLiuZha04} and \citet{AiFanHe07} extended the
FSOA method to construct
$\chi^2$-optimal SSDs. \citet{KouMan05}
constructed some $E(f_{\mathrm{NOD}})$-optimal mixed-level SSDs
by juxtaposing either a saturated two-level orthogonal array and an
$E(f_{\mathrm{NOD}})$-optimal
mixed-level SSD, or two $E(f_{\mathrm{NOD}})$-optimal SSDs.
\citet{Fanetal04N1} and \citet{Tanetal07} presented some
methods for constructing $E(f_{\mathrm{NOD}})$- and $\chi^2$-optimal
mixed-level SSDs, respectively,
from given combinatorial designs. There are many constraints on the
parameters of saturated orthogonal arrays and combinatorial designs
and the construction of most combinatorial designs are unresolved.
Thus, the optimal SSDs obtained by their methods are rather limited.

(c) Yamada et al. (\citeyear{Yametal06}) presented
a method for constructing mixed-level SSDs by juxtaposing two SSDs,
each of which is generated by the
operation ``$\oplus$'' of an initial matrix and a generating matrix. It
can be seen that
their operation ``$\oplus$'' is in fact equivalent to the ``$\oplus
_{\mathcal{A}}$''
in this paper with $\mathcal{A}=\{0, \ldots, q-1\}$, and they only provided
the theoretical justification of the $\chi^2$-optimality for the SSD
with $n=6$.
Recently, \citet{LiuLin09} proposed a method to construct $\chi
^2$-optimal mixed-level SSDs from smaller multi-level SSDs and
transposed orthogonal arrays based on Kronecker sums. It can be easily
confirmed that the result of \citet{LiuLin09} is merely a special
case of our Theorem \ref{sthm2}, by taking $F_1$ as $L_{q_1}=(0, \ldots
, q_1-1)'$ and $D$ as $L_{rq_1}(q_1^{n_2})$. Thus, all their designs
can be constructed by our Theorem \ref{sthm2}.

(d) Using $k$-cyclic generators, \citet{CheLiu08N1} and
\citet{LiuZha09} constructed
some $E(f_{\mathrm{NOD}})$- and $\chi^2$-optimal mixed-level SSDs,
respectively. The
$k$-cyclic generators were obtained via computer searches, when the
values of $k$, the run size and/or the level sizes become larger, the
computer searches tend to be ineffective and impractical.

(e) Recently, \citet{LiuCai09} proposed a new construction
method, called the substitution
method, for $E(f_{\mathrm{NOD}})$-optimal SSDs. It can be seen that all the
$E(f_{\mathrm{NOD}})$-optimal SSDs tabulated in our Tables \ref{tt1} and \ref
{tt3} are different from those tabulated in their Appendices.

Note that the newly proposed methods use small equidistant
designs and difference matrices to generate large designs. Many
difference matrices can be found in Hedayat, Slone and Stufken
(\citeyear{HedSloStu99}), \citet{WuHam00} and from the site
\href{http://support.sas.com/techsup/technote/ts723.html}{http://support.sas.com/techsup/technote/}%
\href{http://support.sas.com/techsup/technote/ts723.html}{ts723.html} maintai\-ned by
Dr. W. F. Kuhfeld of SAS. Equidistant designs can be found
in Ngugen (\citeyear{Ngu96}),
\citet{TanWu97}, \citet{LiuZha00}, Lu et~al. (\citeyear
{Luetal}), \citet{FanLinLiu03}, Fang,
Ge and Liu (\citeyear{FanGeLiu02N1,FanGeLiu02N2,FanGeLiu04}), \citet{LuHuZhe03},
Fang et~al. (\citeyear{Fanetal03,Fanetal04N1,Fanetal04N2}), \citet{AggGup04}, \citet{Esketal04},
\citet{GeoKou06},
Georgiou, Koukouvinos and Mantas (\citeyear{GeoKouMan06}), \citet
{CheLiu08N1}, \citet{LiuCai09} and others.
Difference matrices can also be obtained from orthogonal arrays or
by taking the Kronecker sums of difference matrices, and equidistant
designs also include saturated orthogonal arrays of strength
two.\looseness=1

The appealing feature of our methods is that they can be easily applied
and the resulting designs are $E(f_{\mathrm{NOD}})$- and/or
$\chi^2$-optimal SSDs without fully aliased columns. In particular,
the nonorthogonality between columns of the resulting designs is
well-controlled by the source designs, that is, if the source designs
have little nonorthogonality, the generated design will also have
little nonorthogonality. From these proposed methods, many optimal
SSDs can be constructed in addition to those tabulated in the \hyperref[appm]{Appendix}.

In regard to the statistical data analysis for mixed-level SSDs, it
should be noted that
analyzing the data collected by such SSDs is a very important but
complicated task which has attracted much recent attention. See, for
example, \citet{ZhaZhaLiu07}, \citet{PhoPanXu09} and
\citet{LiZhaZha10}. When there are many more factors than the
number of permitted runs due to expense (e.g., money or time), the
nonorthogonality among factors may be very severe and may prevent the
few active factors to be identified correctly by any existing method.
Therefore, the data analysis for SSDs in general remains an important
and challenging topic for further research. Some recent study on the
analysis of ``high-dimension and low-sample size'' in genetic studies
(e.g., studying 6,000 genes with only 37 observations) may be relevant.

\begin{appendix}\label{appm}\vspace*{2pt}
\section{\texorpdfstring{Proofs}{Appendix A: Proofs}}\label{apdxa}\vspace*{-12pt}
\begin{pf*}{Proof of Theorem \ref{slem0}}
 (a)(i) If $f^1$ and $f^3$ are
fully aliased, that is, $f^1$ can be obtained by permutating the levels
of $f^3$, then there must exist a permutation matrix $Q$ of order
$q$ that satisfies $X_1=X_3Q$, thus $X_1X'_1=X_3QQ'X'_3=X_3X'_3$.

On the other hand, let $V_1$ and $V_3$ be the vector spaces spanned
by the columns of $X_1$ and $X_3$, respectively. If
$X_1X'_1=X_3X'_3$, then $V_1=V_3$, and for any column $x^0$ of
$X_1$, we have
\[
x^0=k_1x^1_3+\cdots+k_qx^q_3  \qquad \mbox{where } x^i_3 \mbox{ is the }
i\mbox{th column of } X_3, i=1, \ldots, q.
\]
Since any two columns in an induced matrix share no element 1 at
any position and each column has
$n_1/q$ ones and $n_1-n_1/q$ zeros, there
must exist only one $k_i\neq0$, that is, $x^0$ is identical to a
column of $X_3$. Then there exists a permutation matrix $Q$ of order
$q$ satisfying $X_1=X_3Q$, thus $f^1$ can be obtained by permutating
the levels of $f^3$, that is, $f^1$ and $f^3$ are fully aliased.

(ii) Note that the induced matrix of $f^1\oplus_{\mathcal{A}}f^2$ is
$[(X_2P_{f_{11}})', \ldots, (X_2P_{f_{n_11}})']'$ and
$P_{f_{i1}}=\sum^{q}_{t=1} x^1_{it}P_{t-1}$, then
\begin{eqnarray*}
\pmatrix{\displaystyle
X_2P_{f_{11}} \cr\displaystyle
\vdots\cr\displaystyle
X_2P_{f_{n_11}}
}
&=&\operatorname{diag}\{\underbrace{X_2, \ldots,
X_2}_{n_1}\}\pmatrix{\displaystyle
\sum^{q}_{t=1} x^1_{1t}P_{t-1} \cr\displaystyle
\vdots\cr\displaystyle
\sum^{q}_{t=1} x^1_{n_1t}P_{t-1}
} \\
&=&(I_{n_1}\otimes X_2)(X_1\otimes I_q)P =(X_1\otimes X_2)P,
\end{eqnarray*}
where $I_n$ is the
identity matrix of order $n$.

(iii) The induced matrices of $f^1\oplus_{\mathcal{A}}f^2$ and
$f^3\oplus_{\mathcal{A}}f^4$ are $(X_1\otimes X_2)P$ and
$(X_3\otimes X_4)P$, respectively. If $f^1\oplus_{\mathcal{A}}f^2$
and $f^3\oplus_{\mathcal{A}}f^4$ are fully aliased, then there
exists a permutation matrix $Q$ of order $q$ such that
\[
(X_1\otimes X_2)P=(X_3\otimes X_4)PQ, \qquad     \mbox{that is, }
\sum^{q}_{t=1}x^1_{st}X_2P_{t-1}=\sum^{q}_{t=1}x^3_{st}X_4P_{t-1}Q
\]
for $s=1, \ldots, n_1.$
For any $s$ and $i$, there is only one nonzero element of
$x^i_{st}$ for $t=1, \ldots, q$ that equals 1. Thus,
$X_2P_{t_1-1}=X_4P_{t_3-1}Q$, and therefore $f^2$ is fully aliased with
$f^4$. Similarly since $f^2\oplus_{\mathcal{A}}f^1$ and
$f^4\oplus_{\mathcal{A}}f^3$ are also fully aliased, it follows
that $f^1$ is fully aliased with $f^3$.

(b)(i) It can be obtained easily from the definition of an induced
matrix.

(ii) From (a)(i), we only need to prove that columns
$q_2(f^1-\frac{q_1-1}{2})\oplus(f^2-\frac{q_2-1}{2})+\frac{q_1q_2-1}{2}$
and
$q_4(f^3-\frac{q_3-1}{2})\oplus(f^4-\frac{q_4-1}{2})+\frac{q_3q_4-1}{2}$
are fully\vspace*{1pt} aliased if and only if $X_1X'_1=X_3X'_3$ and
$X_2X'_2=X_4X'_4$. From (b)(i), the induced matrices of these two
columns are $X_1\otimes X_2$ and $X_3\otimes X_4$, respectively,
thus from (a)(i), they are fully aliased if and only if $(X_1\otimes
X_2)(X_1\otimes X_2)'=(X_3\otimes X_4)(X_3\otimes X_4)'$, that is,
$X_1X'_1\otimes X_2X'_2=X_3X'_3\otimes X_4X'_4$, which means that
$X_1X'_1=aX_3X'_3$ and $X_2X'_2=1/aX_4X'_4$ for some $a\neq0$. Since
the elements in $X_iX'_i$ are all ones and zeros for $i=1, \ldots,
4,$ then $a=1$, that is, $X_1X'_1=X_3X'_3$ and $X_2X'_2=X_4X'_4$.

(iii) The induced matrices of columns
$q_2(f^1-\frac{q_1-1}{2})\oplus(f^2-\frac{q_2-1}{2})+\frac{q_1q_2-1}{2}$
and $f^3\oplus_{\mathcal{A}}f^4$ are $X_1\otimes X_2 \mbox{ and }
[(X_4P_{f_{13}})', \ldots, (X_4P_{f_{n_13}})']',$ respectively. If
these two columns are fully aliased, then
$ v_{ij}X_2X'_2=X_4P_{f_{i3}}P'_{f_{j3}}X'_4,$ where $v_{ij}$
is the $(i,j)$th entry of $X_1X'_1, i, j=1,\ldots, n_1.$
Note that $v_{ij}$ can be zero, and hence $v_{ij}X_2X'_2$ can be a
zero matrix which contradicts the fact that
$X_4P_{f_{i3}}P'_{f_{j3}}X'_4$ cannot be a zero matrix in any case.
\end{pf*}

\begin{pf*}{Proof of Theorem \ref{sthm1}}
 (a) Consider the $i$th and $j$th
rows of $F\oplus_{\mathcal{A}}D'$,
$(f_{i_1}\oplus_{\mathcal{A}}d_{i_2})'$ and
$(f_{j_1}\oplus_{\mathcal{A}}d_{j_2})'$, where $i=(i_1-1)c+i_2,
j=(j_1-1)c+j_2, i_1, j_1=1, \ldots, n, i_2, j_2=1, \ldots, c,$ and
$i\neq j$, $f_k$ and $d_k$ are the $k$th rows of $F$ and $D'$,
respectively. Then the coincidence number between
$(f_{i_1}\oplus_{\mathcal{A}}d_{i_2})'$ and
$(f_{j_1}\oplus_{\mathcal{A}}d_{j_2})'$ equals the number of zeros
in $(f_{i_1}-f_{j_1})\oplus_{\mathcal{A}}(d_{i_2}-d_{j_2})$.

(i) Suppose $i_1=j_1, i_2\neq j_2$, then $f_{i_1}=f_{j_1}$ and
$d_{i_2}\neq d_{j_2}$. From the definition of difference matrix,
each element in $\mathcal{A}$ occurs $r$ times in $d_{i_2}-d_{j_2}$.
Therefore, $(f_{i_1}-f_{j_1})\oplus_{\mathcal{A}}(d_{i_2}-d_{j_2})
=0_m\oplus_{\mathcal{A}}(d_{i_2}-d_{j_2})$ and there are $mr$ zeros
in $0_m\oplus_{\mathcal{A}}(d_{i_2}-d_{j_2})$, where $0_m$ denotes
the $m\times1$ column vector with all elements zero, that is, the
coincidence number between
$(f_{i_1}\oplus_{\mathcal{A}}d_{i_2})'$ and
$(f_{j_1}\oplus_{\mathcal{A}}d_{j_2})'$ is $mr$.

(ii) If $i_1\neq j_1, i_2\neq j_2$, similar to (i), it can also
be easily seen that there are $mr$ zeros in
$(f_{i_1}-f_{j_1})\oplus_{\mathcal{A}}(d_{i_2}-d_{j_2})$, that is,
the coincidence number between
$(f_{i_1}\oplus_{\mathcal{A}}d_{i_2})'$ and
$(f_{j_1}\oplus_{\mathcal{A}}d_{j_2})'$ is $mr$.

(iii) If $i_1\neq j_1, i_2=j_2$, that is, $f_{i_1}\neq f_{j_1}$,
$d_{i_2}=d_{j_2},$ then
$(f_{i_1}-f_{j_1})\oplus_{\mathcal{A}}(d_{i_2}-d_{j_2})=
(f_{i_1}-f_{j_1})\oplus_{\mathcal{A}}0_{rq}$, and there are
$\lambda rq$ zeros in
$(f_{i_1}-f_{j_1})\oplus_{\mathcal{A}}0_{rq}$, that is, the
coincidence number between
$(f_{i_1}\oplus_{\mathcal{A}}d_{i_2})'$ and
$(f_{j_1}\oplus_{\mathcal{A}}d_{j_2})'$ is $\lambda rq$.

(b) $F\oplus_{\mathcal{A}}D'$ can be obtained from
$D'\oplus_{\mathcal{A}}F$ through row and column permutations. Thus,
if $D'\oplus_{\mathcal{A}}F$ has no fully aliased columns, neither
does $F\oplus_{\mathcal{A}}D'$. Let $d^1\oplus_{\mathcal{A}}f^1$ and
$d^2\oplus_{\mathcal{A}}f^2$ be two different columns of
$D'\oplus_{\mathcal{A}}F$, where $d^i=(d_{1i}, \ldots, d_{ci})'
\mbox{ and } f^i=(f_{1i}, \ldots, f_{ni})' \mbox{ for } i=1 \mbox{
and } 2$ are columns of $D'$ and $F$, respectively. Since $D$ is a
normalized difference matrix, $d_{1i}=0$ for $i=1 \mbox{ and } 2.$
Let $X_1$ and $X_2$ be the induced matrices of $f^1$ and $f^2$,
respectively. Then the induced matrices of
$d^1\oplus_{\mathcal{A}}f^1$ and $d^2\oplus_{\mathcal{A}}f^2$ are
\[
 \pmatrix{\displaystyle
X_1 \cr\displaystyle
X_1P_{d_{21}} \cr\displaystyle
\vdots\cr\displaystyle
X_1P_{d_{c1}}
}
 \quad \mbox{and} \quad
 \pmatrix{\displaystyle
X_2 \cr\displaystyle
X_2P_{d_{22}} \cr\displaystyle
\vdots\cr\displaystyle
X_2P_{d_{c2}} \cr\displaystyle
} \qquad   \mbox{respectively}.
\]
Suppose
$d^1\oplus_{\mathcal{A}}f^1$ and $d^2\oplus_{\mathcal{A}}f^2$ are
fully aliased. Then from Theorem \ref{slem0},
%
\begin{equation}\label{1}
X_1P_{d_{1i}}P'_{d_{j1}}X'_1=X_2P_{d_{i2}}P'_{d_{j2}}X'_2,\qquad  i,
j=1,\ldots, c.
\end{equation}
Noting that $P_{d_{11}}=P_{d_{12}}=I_q$, we can obtain the following
equations by taking $i=1$ in \eqref{1}:
%
\begin{eqnarray} \label{2}
X_1X'_1&=&X_2X'_2, \\
\label{21}
X_1P'_{d_{j1}}X'_1&=&X_2P'_{d_{j2}}X'_2, \qquad j=2,\ldots, c.
\end{eqnarray}
Since $F$ has no fully aliased columns, from equation \eqref{2}, we
know that $f^1$ and $f^2$ must be the same column of $F$, thus
$X_1=X_2$, and $X_1P'_{d_{j1}}X'_1=X_1P'_{d_{j2}}X'_1$ for
$j=2,\ldots, c.$ Also, since $X_1$ is a column full rank matrix,\vspace*{-1pt}
we have $P'_{d_{j1}}=P'_{d_{j2}}$, and thus $d_{j1}=d_{j2}, \mbox{
for } j=1,\ldots, c$, that is, $d^1=d^2$. So $d^1$ and $d^2$ must be
the same row of $D$ since $D$ has no identical rows. Therefore,
$d^1\oplus_{\mathcal{A}}f^1$ and $d^2\oplus_{\mathcal{A}}f^2$ are the
same column of $D'\oplus_{\mathcal{A}} F,$ which contradicts the
fact that they are two different columns of $D'\oplus_{\mathcal{A}}
F.$ Hence, $D'\oplus_{\mathcal{A}} F$ as well as
$F\oplus_{\mathcal{A}} D'$ have no fully aliased columns.
\end{pf*}

\begin{pf*}{Proof of Theorem \ref{sthm2}} We only prove that there are
no fully aliased columns between $F_1\oplus_{\mathcal{A}_1}D'$ and
$0_{n_1}\oplus F_2$. (The others can be proved easily.) Suppose
$f^1\oplus_{\mathcal{A}_1}d^1$ and $0_{n_1}\oplus f^2$ are columns
of $F_1\oplus_{\mathcal{A}_1}D'$ and $0_{n_1}\oplus F_2$,
respectively, where $f^1$, $f^2$ and $d^1=(0, d_{21}, \ldots,
d_{n_21})'$ are columns of $F_1$, $F_2$ and $D'$, respectively. Let
$X$ and $Y$ be the induced matrices of $f^1$ and $f^2$,
respectively. Then the induced matrices of
$d^1\oplus_{\mathcal{A}_1}f^1$ and $f^2\oplus0_{n_1}$ are $[X',
(XP_{d_{21}})', \ldots, (XP_{d_{n_21}})'] \mbox{ and } Y\otimes
1_{n_1},$ respectively. From the definition of an induced matrix, it
is easy to see that $XX'\neq y_01_{n_1}1'_{n_1}$,\vspace*{1pt} where $y_0$ is the
$(1, 1)$th entry of $YY'$ and $1_{n_1}$ denotes the $n_1\times1$
vector with all elements unity. Thus, from Theorem \ref{slem0},
$d^1\oplus_{\mathcal{A}_1}f^1$ and $f^2\oplus0_{n_1}$ are not fully
aliased. Therefore, $f^1\oplus_{\mathcal{A}_1}d^1$ and $0_{n_1}\oplus
f^2$ are not fully aliased.
\end{pf*}

The following lemma will be used in the proof of Theorem
\ref{sthm6}.\vadjust{\goodbreak}

\begin{lema}[{[\citet{FanLinLiu03}]}]\label{slem5}
Suppose $f^j$ is the $j$th column of an $F(n, q_1\cdots q_m)$ with
induced matrix $X_j$, $j=1, \ldots, m$. Then
\[
f_{\mathrm{NOD}}(f^i, f^j)=\operatorname{tr}(X'_iX_jX'_jX_i)-\frac{n^2}{q_iq_j}.
\]
\end{lema}

\begin{pf*}{Proof of Theorem \ref{sthm6}}
(a) From Theorem \ref{slem0} and
Lemma \ref{slem5}, the induced matrices of $h_1$ and $h_2$ are
$X_1\otimes X_2$ and $X_3\otimes X_4$, respectively. Then we have
\begin{eqnarray*}
f_{\mathrm{NOD}}(h_1,
h_2)&=&\operatorname{tr}[(X'_1X_3X'_3X_1)\otimes(X'_2X_4X'_4X_2)]-\frac
{n^2_1n^2_2}{\prod^{4}_{i=1}q_i} \\
&=&\operatorname{tr}(X'_1X_3X'_3X_1)\operatorname{tr}(X'_2X_4X'_4X_2)-\frac
{n^2_1n^2_2}{\prod^{4}_{i=1}q_i} \\
&=&f_{\mathrm{NOD}}(f^1,f^3)f_{\mathrm{NOD}}(f^2, f^4)\\
&&{}+\frac{n^2_2}{q_2q_4}f_{\mathrm{NOD}}(f^1,f^3)+\frac
{n^2_1}{q_1q_3}f_{\mathrm{NOD}}(f^2, f^4).
\end{eqnarray*}

(b) The induced matrices of $f^1\oplus_{\mathcal{A}_1}f^2$ and
$f^3\oplus_{\mathcal{A}_3}f^4$ are $(X_1\otimes X_2)P$ and
$(X_3\otimes X_4)Q$, respectively, where $P=(P'_0, \ldots,
P'_{q_1-1})', Q=(Q'_0, \ldots, Q'_{q_3-1})'$, $P_i$ and $Q_j$ are
permutation matrices defined by $(0, \ldots,
q_1-1)P'_i=i+_{{\mathcal{A}_1}}(0, \ldots, q_1-1)$ and $(0, \ldots,
q_3-1)Q'_j=j+_{{\mathcal{A}_3}}(0, \ldots, q_3-1),$ respectively,
$i=0, \ldots, q_1-1, j=0, \ldots, q_3-1$. Let $T=P'(X'_1X_3\otimes
X'_2X_4)Q-\frac{n_1n_2}{q_1q_3}{1}_{q_1}{1}'_{q_3}$. Then from Lemma
\ref{slem5}$, f_{\mathrm{NOD}}(f^1\oplus_{\mathcal{A}_1}f^2,
f^3\oplus_{\mathcal{A}_3}f^4)$ equals the sum of squares of the
elements of $T$. Let $W=(w_{ij})=X'_1X_3,
B=(b_{ij})=X'_2X_4-n_2/(q_1q_3){1}_{q_1}{1}'_{q_3}$, and note that
$\sum^{q_3}_{j=1}\sum^{q_1}_{i=1}w_{ij}=n_1,
\sum^{q_3}_{j=1}\sum^{q_1}_{i=1}b_{ij}=0$. Then
\[
T=\sum^{q_3}_{j=1}\sum^{q_1}_{i=1}w_{ij}P'_{i-1}BQ_{j-1}
\]
and the $(s, t)$th entry of $T$ can be expressed as\vspace*{1pt}
$\sum^{q_3}_{j=1}\sum^{q_1}_{i=1}w_{ij}b_{s_i t_j}$, where $(s_1,
\ldots, s_{q_1})$ and $(t_1, \ldots, t_{q_3})$ are some permutations
of $(1, \ldots, q_1)$ and $(1, \ldots,\break q_3)$, respectively. Then
\[
 \Biggl(\sum^{q_3}_{j=1}\sum^{q_1}_{i=1}w_{ij}b_{s_i t_j} \Biggr)^2\leq
 \Biggl(\sum^{q_3}_{j=1}\sum^{q_1}_{i=1}w^2_{ij} \Biggr)\Biggl (\sum
^{q_3}_{j=1}\sum^{q_1}_{i=1}b^2_{ij} \Biggr),
\]
and thus
\begin{eqnarray*}
f_{\mathrm{NOD}}(f^1\oplus_{\mathcal{A}_1}f^2,
f^3\oplus_{\mathcal{A}_3}f^4)&\leq&
q_1q_3 \Biggl(\sum^{q_3}_{j=1}\sum^{q_1}_{i=1}w^2_{ij} \Biggr)\Biggl (\sum
^{q_3}_{j=1}\sum^{q_1}_{i=1}b^2_{ij} \Biggr) \\
&=&q_1q_3 \biggl[f_{\mathrm{NOD}}(f^1,
f^3)+\frac{n^2_1}{q_1q_3} \biggr]f_{\mathrm{NOD}}(f^2, f^4),
\end{eqnarray*}
where the equality holds if and only if there exist $c_1$ and
$c_2$ with $|c_1|+|c_2|>0$ such that $c_1w_{ij}=c_2b_{s_i t_j}$ for
$i=1, \ldots, q_1$ and $j=1, \ldots, q_3$. This means that
$c_1\sum^{q_3}_{j=1}\sum^{q_1}_{i=1}w_{ij}=c_2\sum^{q_3}_{j=1}\sum
^{q_1}_{i=1}b_{s_i
t_j}=0$, and so $c_1=0, c_2\neq0$ and $b_{ij}=0$ for $i=1, \ldots,
q_1$ and $j=1, \ldots, q_3.$ Thus, $f^2$ is orthogonal to $f^4$.

On the other hand, if $f^2$ is orthogonal to $f^4$, $f_{\mathrm{NOD}}(f^2,
f^4)=0$ and
\begin{eqnarray*}
0&\leq& f_{\mathrm{NOD}}(f^1\oplus_{\mathcal{A}_1}f^2,
f^3\oplus_{\mathcal{A}_3}f^4)\\
&\leq& q_1q_3\biggl [f_{\mathrm{NOD}}(f^1,
f^3)+\frac{n^2_1}{q_1q_3} \biggr]f_{\mathrm{NOD}}(f^2, f^4)=0,
\end{eqnarray*}
then the
equality holds.

Similarly, we can obtain that
\begin{eqnarray*}
f_{\mathrm{NOD}}(f^1\oplus_{\mathcal{A}_1}f^2,
f^3\oplus_{\mathcal{A}_3}f^4)
&=&f_{\mathrm{NOD}}(f^2\oplus_{\mathcal{A}_1}f^1,f^4\oplus_{\mathcal
{A}_3}f^3)\\
&\leq&
q_2q_4 \biggl[f_{\mathrm{NOD}}(f^2,f^4)+\frac{n^2_2}{q_2q_4} \biggr]f_{\mathrm{NOD}}(f^1,
f^3),
\end{eqnarray*}
and the equality holds if and only if $f^1$ is orthogonal to $f^3$.
Hence, we have the
assertion.

(c) The induced matrices of $f^1\oplus_{\mathcal{A}_1}f^2$ and
$h_2$ are $(X_1\otimes X_2)P$ and $X_3\otimes X_4$, respectively.
Let $K=P'(X'_1X_3\otimes
X'_2X_4)-\frac{n_1n_2}{q_1q_3q_4}{1}_{q_1}{1}'_{q_3q_4},
G=(g_{ij})=X'_2X_4-\frac{n_2}{q_1q_4}{1}_{q_1}{1}'_{q_4}$ and
$W=(w_{ij})=X'_1X_3$, and note that $\sum^{q_1}_{i=1}w_{ij}=n_1/q_3,
j=1, \ldots, q_3$. Then
$K=(A_1, \ldots, A_{q_3}), \mbox{ where } A_j=\sum^{q_1}_{i=1}w_{ij}P'_{j-1}G.$
Note that $f_{\mathrm{NOD}}(f^1\oplus_{\mathcal{A}_1}f^2, h_2)$ is equal
to the sum of squares of the elements of $K$, the $(s, t)$th entry
of $A_j$ is $\sum^{q_1}_{i=1}w_{ij}g_{s_i t}$ and\vspace*{1pt}
$(\sum^{q_1}_{i=1}w_{ij}g_{s_i t})^2\leq
\sum^{q_1}_{i=1}w^2_{ij}\sum^{q_1}_{s=1}g^2_{st}$, where $(s_1,
\ldots, s_{q_1})$ is a permutation of $(1, \ldots, q_1)$. Then
similar to the proof in (b), we get
\begin{eqnarray*}
f_{\mathrm{NOD}}(f^1\oplus_{\mathcal{A}_1}f^2, h_2)&\leq& \sum
^{q_3}_{j=1}\sum^{q_4}_{t=1}
\sum^{q_1}_{s=1} \Biggl(\sum^{q_1}_{i=1}w^2_{ij}\sum
^{q_1}_{k=1}g^2_{kt} \Biggr)\\
&=&q_1\sum^{q_3}_{j=1}\sum^{q_1}_{i=1}w^2_{ij}\sum^{q_4}_{t=1}\sum
^{q_1}_{k=1}g^2_{kt} \\
&=&q_1\biggl [f_{\mathrm{NOD}}(f^1, f^3)+\frac{n^2_1}{q_1q_3}
\biggr]f_{\mathrm{NOD}}(f^2, f^4),
\end{eqnarray*}
where the equality holds if and only if $f^2$ is orthogonal to
$f^4$, and
\begin{eqnarray*}
f_{\mathrm{NOD}}(f^1\oplus_{\mathcal{A}_1}f^2, h_2)\leq q_1
\biggl[f_{\mathrm{NOD}}(f^2, f^4)+\frac{n^2_2}{q_1q_4} \biggr]f_{\mathrm{NOD}}(f^1, f^3),
\end{eqnarray*}
where the equality holds if and only if $f^1$ is orthogonal to
$f^3$. Thus, we complete the proof of (c).
\end{pf*}

\section{\texorpdfstring{Some selected optimal supersaturated designs}%
{Appendix B: Some selected optimal supersaturated designs}}
\label{apdxb}

\begin{table}[b]
\vspace{-20pt}
\caption{Equidistant designs used in Tables \protect\ref{tt1}--\protect\ref{tt2}}
\tabcolsep=0pt
\begin{tabular*}{\textwidth}{@{\extracolsep{\fill}}lccl@{}}
\hline
$\bolds{n}$&$\bolds{m}$&$\bolds{q}$& \multicolumn{1}{c@{}}{\textbf{Source design}} \\
 \hline
\hphantom{2}4 & 3 &  2 & Orthogonal array \\
\hphantom{2}8 & 7 &  2 & Orthogonal array \\
12 & 11 &  2 & Orthogonal array \\
16 & 15 &  2 & Orthogonal array \\
16 & 5 &  4 & Orthogonal array \\
20 & 19 &  2 & Orthogonal array \\
24 & 23 &  2 & Orthogonal array \\
25 & 6 &  5 & Orthogonal array \\
\hphantom{2}6 & 10 &  2 & \citet{LiuZha00} \\
\hphantom{2}6 & 5 & 3 & \citet{FanGeLiu04} \\
\hphantom{2}6 & $5k$  $(k=2, 3)$ & 3 & \citet{GeoKou06} \\
\hphantom{2}8 & $7k$  $(k=2, \ldots, 5)$ & 2 & \citet{LiuZha00} \\
\hphantom{2}8 & $7k$   $(k=1, 2)$ & 4 & \citet{FanGeLiu02N1} \\
\hphantom{2}8 & $7k$   $(k=3, \ldots, 6)$ & 4 & \citet{GeoKou06} \\
\hphantom{2}9 & $4k$   $(k=1, \ldots, 7)$ & 3 & \citet{FanGeLiu04} \\
\hphantom{2}9 & $4k$   $(k=8, 10, 12)$ & 3 & \citet{GeoKou06} \\
10 & $18k$   $(k=1, 2, 3)$ & 2 & \citet{LiuZha00} \\
10 & $9$   & 5 & \citet{FanGeLiu02N2} \\
10 & $9k$   $(k=2, 3, 4)$ & 5 & \citet{GeoKou06} \\
12 & $11k$   $(k=2, \ldots, 12)$ & 2 & \citet{LiuZha00} \\
12 & $11$ &   3 & \citet{LuHuZhe03} \\
12 & $11k$   $ (k=2, \ldots, 5)$ & 3 & \citet{GeoKou06} \\
12 & $11$ &  6 & \citet{LuHuZhe03} \\
12 & $11k$   $(k=2, 3)$ & 6 & \citet{GeoKou06} \\
14 & $13k$   $(k=1, 2)$ & 7 & \citet{Fanetal03} \\
15 & 28 &  3 & \citet{GeoKou06} \\
15 & $7k$   $(k=1, \ldots, 13)$ & 5 & \citet{FanGeLiu04} \\
16 & $15k$   $(k=2, \ldots, 6)$ & 2 & \citet{LiuZha00} \\
16 & $15k$   $(k=7, 8, 9)$ & 2 & \citet{Esketal04} \\
16 & $5k$   $(k=2, \ldots, 6)$ & 4 & \citet{Fanetal03} \\
16 & $5k$   $(k=7, \ldots, 16)$ & 4 & Georgiou, Koukouvinos and Mantas
(\citeyear{GeoKouMan06}) \\
18 & $34k$   $(k=1, 2, 3)$ & 2 & \citet{LiuZha00} \\
18 & $17k$   $(k=1, 2)$ & 3 & \citet{Fanetal03}\\
18 & $17$ & 6 & \citet{LuHuZhe03} \\
18 & $34$ & 6 & \citet{GeoKou06} \\
20 & $19k$  $(k=2, 3)$ & 2 & \citet{LiuZha00} \\
20 & 19 & 4 & Lu et~al. (\citeyear{Luetal}) \\
20 & 19 & 5 & \citet{Luetal} \\
22 & 42 & 2 & \citet{LiuZha00} \\
24 & 46 & 2 & \citet{LiuZha00} \\
24 & 23 & 4 & \citet{Luetal} \\
24 & 23 & 6 & \citet{LuHuZhe03} \\
25 & $6k$  $(k=2, \ldots, 25)$ & 5 & Georgiou, Koukouvinos and Mantas
(\citeyear{GeoKouMan06}) \\
\hline
\end{tabular*}
\vspace{-15pt}
\end{table}

\begin{table}
\tabcolsep=0pt
\caption{\vspace*{-1pt}Some selected $E(f_{\mathrm{NOD}})$-optimal SSDs constructed by
Theorem \protect\ref{sthm2}}\vspace*{-1pt}
\label{tt1}
{\fontsize{8.8pt}{10.8pt}\selectfont{\begin{tabular*}{\textwidth}{@{\extracolsep{\fill}}ld{2.0}cd{2.0}ccccc}
\hline\\[-12pt]
$\bolds{n_1}$&\multicolumn{1}{c}{$\bolds{m_1}$}&$\bolds{q_1}$& \multicolumn{1}{c}{$\bolds{n_2}$}
&$\bolds{m_2}$&$\bolds{q_2}$& $\bolds{r}$ & \textbf{Final resulting
SSD}\tabnoteref[\dag]{tab1}& $\bolds{\lambda}$\\\hline\\[-12pt]
\hphantom{1}4 & 3 & 2 & 6 & $5k$ & 3 & $4k$ & $F(24, 2^{24k}3^{5k})$ & $13k,\ k=1,
2, 3$ \\
\hphantom{1}4 & 3 & 2 & 8 & $7k$ & 4 & $6k$ & $F(32, 2^{36k}4^{7k})$ & $19k,\ k=1,
\ldots, 6$ \\
\hphantom{1}4 & 3 & 2 & 9 & $4k$ & 3 & $3k$ & $F(36, 2^{18k}3^{4k})$ & $10k,\ k=2t,
t=1\ldots, 6$ \\
\hphantom{1}6 & 5 & 3 & 6 & 10 &2 & 3 & $F(36, 2^{10}3^{45})$ & 19 \\
\hphantom{1}4 & 3 & 2 & 10 & $9k$ & 5 & $8k$ & $F(40, 2^{48k}5^{9k})$ & $25k,\ k=1,
\ldots, 4$ \\
\hphantom{1}4 & 3 & 2 & 12 & $11k$ & 3 & $8k$ & $F(48, 2^{48k}3^{11k})$ & $27k,\ k=1, \ldots, 5$ \\
\hphantom{1}6 & 5 & 3 & 8 & $7k$ & 2 & $2k$ & $F(48, 2^{7k}3^{30k})$ & $13k,\ k=2,
\ldots, 5$ \\
\hphantom{1}{6} & 10 & 2 & 8 & $7k$ & 4 & $3k$ & $F(48, 2^{60k}4^{7k})$ & $31k,\ k=2, 4, 6$ \\
\hphantom{1}4 & 3 & 2 & 12 & $11k$ & 6 & $10k$ & $F(48, 2^{60k}6^{11k})$ & $31k,\ k=1, 2, 3$ \\
\hphantom{1}6 & 5 & 3 & 8 & $7k$ & 4 & $3k$ & $F(48, 3^{45k}4^{7k})$ & $16k,\ k=1,
\ldots, 6$ \\
\hphantom{1}6 & 10 & 3 & 8 & $7k$ & 2 & $k$ & $F(48, 2^{7k}3^{30k})$ & $13k,\ k=3,
4, 5$ \\
\hphantom{1}6 & 10 & 3 & 8 & $14k$ & 4 & $3k$ & $F(48, 3^{90k}4^{14k})$ & $32k,\ k=1, 2, 3$ \\
\hphantom{1}{6} & 10 & 2 & 9 & $16$ & 3 & $6$ & $F(54, 2^{120}3^{16})$ & 64 \\
\hphantom{1}4 & 3 & 2 & 14 & $13k$ & 7 & $12k$ & $F(56, 2^{72k}7^{13k})$ & $37k,\ k=1, 2$ \\
\hphantom{1}4 & 3 & 2 & 15 & $28$ & 3 & $20$ & $F(60, 2^{120}3^{28})$ & 68 \\
\hphantom{1}{6} & 5 & 3 & 10 & $18k$ & 2 & $5k$ & $F(60, 2^{18k}3^{75k})$ & $33k,\ k=2, 3$ \\
\hphantom{1}4 & 3 & 2 & 15 & $7k$ & 5 & $6k$ & $F(60, 2^{36k}5^{7k})$ & $19k,\ k=2,
\ldots, 13$ \\
\hphantom{1}{6} & 10 & 2 & 10 & $9k$ & 5 & $4k$ & $F(60, 2^{80k}5^{9k})$ & $41k,\ k=2, 3, 4$ \\
\hphantom{1}6 & 5 & 3 & 10 & $9k$ &5 & $4k$ & $F(60, 3^{60k}5^{9k})$ & $21k,\ k=1,
\ldots, 4$ \\
{10} & 9 & 5 & 6 & $5k$ & 3 & $k$ & $F(60, 3^{5k}5^{45k})$ & $10k,\ k=2,
3$ \\
\hphantom{1}6 & 10 & 3 & 10 & 18 & 5 & 4 & $F(60, 3^{120}5^{18})$ &42 \\
\hphantom{1}4 & 3 & 2 & 16 & $5k$ & 4 & $4k$ & $F(64, 2^{24k}4^{5k})$ & $13k,\ k=2,
\ldots, 16$ \\
\hphantom{1}4 & 3 & 2 & 18 & $17k$ & 3 & $12k$ & $F(72, 2^{72k}3^{17k})$ & $41k,\ k=1, 2$ \\
\hphantom{1}{6} & 10 & 2 & 12 & $11k$ & 3 & $4k$ & $F(72, 2^{80k}3^{11k})$ & $43k,\ k=2, \ldots, 5$ \\
\hphantom{1}4 & 3 & 2 & 18 & $34$ & 6 & $30$ & $F(72, 2^{180}6^{34})$ & $94$ \\
\hphantom{1}{6} & 10 & 2 & 12 & $22$ & 6 & $10$ & $F(72, 2^{200}6^{22})$ & $102$ \\
\hphantom{1}{9} & 4 & 3 & 8 & $7k$ & 4 & $6k$ & $F(72, 3^{72k}4^{7k})$ & $25k,\ k=1,
\ldots, 6$ \\
\hphantom{1}{6} & 5 & 3 & 12 & $11k$ & 6 & $5k$ & $F(72, 3^{75k}6^{11k})$ & $26k,\ k=2, 3$ \\
{10} & 18 & 2 & 6 & $5k$ & 3 & $2k$ & $F(80, 2^{72k}3^{5k})$ & $37k,\ k=2, 3$ \\
{10} & 18 & 2 & 8 & $7k$ & 4 & $3k$ & $F(80, 2^{108k}4^{7k})$ & $55k,\ k=2, 4, 6$ \\
\hphantom{1}4 & 3 & 2 & 20 & 19 & 5 & 16 & $F(80, 2^{96}5^{19})$ & 51 \\
{10} & 9 & 5 & 8 & $14k$ & 4 & $3k$ & $F(80, 4^{14k}5^{135k})$ & $29k,\ k=1, 2, 3$ \\
\hphantom{1}{6} & 10 & 2 & 14 & 26 & 7 & 12 & $F(84, 2^{240}7^{26})$ & 122 \\
\hphantom{1}6 & 5 & 3 & 14 & $26$ &7 & $12$ & $F(84, 3^{180}7^{26})$ & $62$ \\
\hphantom{1}{6} & 10 & 2 & 15 & $7k$ & 5 & $3k$ & $F(90, 2^{60k}5^{7k})$ & $31k,\ k=2t, t=2, \ldots, 6$ \\
\hphantom{1}6 & 5 & 3 & 15 & $7k$ &5 & $3k$ & $F(90, 3^{45k}5^{7k})$ & $ 16k,\ k=2\ldots, 13$ \\
\hphantom{1}9 & 4 & 3 & 10& $9k$ & 5 & $8k$ & $F(90, 3^{96k}5^{9k})$ & $33k,\ k=1,
\ldots, 4$ \\
\hphantom{1}{6} & 5 & 3 & 16 & $15k$ & 2 & $4k$ & $F(96, 2^{15k}3^{60k})$ & $27k,\ k=2, \ldots, 9$ \\
\hphantom{1}4 & 3 & 2 & 24 & 23 & 4 & 18 & $F(96, 2^{108}4^{23})$ & 59\\
\hphantom{1}{6} & 10 & 2 & 16 & $5k$ & 4 & $2k$ & $F(96, 2^{40k}4^{5k})$ & $21k,\ k=4, \ldots, 16$ \\
\hphantom{1}4 & 3 & 2 & 24 & 23 & 6 & 20 & $F(96, 2^{120}6^{23})$ & 63 \\
\hphantom{1}6 & 5 & 3 & 16 & $5k$ &4 & $2k$ & $F(96, 3^{30k}4^{5k})$ & $11k,\ k=3,
\ldots, 16$ \\
\hphantom{1}4 & 3 & 2 & 25 & $6k$ & 5 & $5k$ & $F(100, 2^{30k}5^{6k})$ & $16k,\ k=2t, t=2, \ldots, 12$ \\[-1pt]
\hline\\[-13pt]
\end{tabular*}}}
\tabnotetext[\dag]{tab1}{$F(n_1n_2, q^{rm_1q_1}_1q^{m_2})$.}
\legend{$\lambda$ is the constant coincidence number
of the final resulting SSD.}
\end{table}

\begin{table}
\tabcolsep=0pt
\caption{Some selected $\chi^2$-optimal SSDs constructed by Theorem
\protect\ref{sthm2}}
\begin{tabular*}{\textwidth}{@{\extracolsep{\fill}}ld{2.0}cd{2.0}ccccc@{}}
\hline
$\bolds{n_1}$&\multicolumn{1}{c}{$\bolds{m_1}$}&$\bolds{q_1}$& \multicolumn{1}{c}{$\bolds{n_2}$}&$\bolds{m_2}$
&$\bolds{q_2}$& $\bolds{r}$ & \textbf{Final resulting
SSD\tabnoteref[\dag]{tab2}} & $\bolds{\omega}$\\
\hline
\hphantom{1}4 & 3 & 2 & 6 & $5k$ & 3 & $6k$ & $F(24, 2^{36k}3^{5k})$ & $39k,\ k=1,
2, 3$ \\
\hphantom{1}4 & 3 & 2 & 8 & $7k$ & 4 & $12k$ & $F(32, 2^{72k}4^{7k})$ & $76k,\ k=1,
\ldots, 6$ \\
\hphantom{1}4 & 3 & 2 & 9 & $16k$ & 3 & $18k$ & $F(36, 2^{108k}3^{16k}$) & $120k,\ k=1, 2, 3$ \\
\hphantom{1}6 & 10 & 2 & 6 & 10 & 3 & 6 & $F(36, 2^{120}3^{10})$ & 126 \\
\hphantom{1}6 & 5 & 3 & 6 & 10 & 2 & 2 & $F(36, 2^{10}3^{30})$ & 38 \\
\hphantom{1}4 & 3 & 2 & 10 & $9k$ & 5 & $20k$ & $F(40, 2^{120k}5^{9k})$ & $125k,\ k=1, \ldots, 4$ \\
\hphantom{1}4 & 3 & 2 & 12 & $11k$ & 3 & $12k$ & $F(48, 2^{72k}3^{11k})$ & $81k,\ k=1, \ldots, 5$ \\
\hphantom{1}8 & 14 & 2 & 6 & $10$ & 3 & $6$ & $F(48, 2^{168}3^{10})$ & $174$\\
\hphantom{1}4 & 3 & 2 & 12 & $11k$ & 6 & $30k$ & $F(48, 2^{180k}6^{11k})$ & $186k,\ k=1, 2, 3$ \\
\hphantom{1}6 & 10 & 2 & 8 & $7k$ & 4 & $6k$ & $F(48, 2^{120k}4^{7k})$ & $124k,\ k=1, \ldots, 6$ \\
\hphantom{1}6 & 5 & 3 & 8 & $7k$ & 4 & $4k$ & $F(48, 3^{60k}4^{7k})$ & $64k,\ k=1,
\ldots, 6$ \\
\hphantom{1}6 & 10 & 3 & 8 & $7k$ & 4 & $2k$ & $F(48, 3^{60k}4^{7k})$ & $64k,\ k=2,
\ldots, 6$ \\
\hphantom{1}6 & 15 & 3 & 8 & 42 & 4 & 8 & $F(48, 3^{360}4^{42})$ & 384 \\
\hphantom{1}8 & 7 & 2 & 6 & $5k$ & 3 & $6k$ & $F(48, 2^{84k}3^{5k})$ & $87k,\ k=1,
2, 3$ \\
\hphantom{1}8 & 7 & 4 & 6 & $5k$ & 3 & $k$ & $F(48, 4^{28k}3^{5k})$ & $31k,\ k=2,
3$\\
\hphantom{1}4 & 3 & 2 & 14 & 13 & 7 & 42 & $F(56, 2^{252}7^{13})$ & 259 \\
\hphantom{1}4 & 3 & 2 & 15 & 28 & 3 & 30 & $F(60, 2^{180}3^{28})$ & 204 \\
\hphantom{1}4 & 3 & 2 & 15 & $14$ & 5 & $30$ & $F(60, 2^{180}5^{14})$ & $190$ \\
\hphantom{1}6 & 10 & 2 & 10 & $9k$ & 5 & $10k$ & $F(60, 2^{200k}5^{9k})$ & $205k,\ k=1, \ldots, 4$ \\
10 & 18 & 2 & 6 & 10 & 3 & 6 & $F(60, 2^{216}3^{10})$ & 222\\
\hphantom{1}4 & 3 & 2 & 16 & $5k$ & 4 & $8k$ & $F(64, 2^{48k}4^{5k})$ & $52k,\ k=2,
\ldots, 16$\\
\hphantom{1}8 & 14& 2 & 8 & $7k$ & 4 & $6k$ & $F(64, 2^{168k}4^{7k})$ & $172k,\ k=1, \ldots, 6$ \\
\hphantom{1}4 & 3 & 2 & 18 & $17k$ & 3 & $18k$ & $F(72, 2^{108k}3^{17k})$ & $123k,\ k=1, 2$\\
\hphantom{1}8 & 7 & 2 & 9 & $8k$ & 3 & $9k$ & $F(72, 2^{126k}3^{8k})$ & $132k,\ k=2, 4, 6$ \\
\hphantom{1}9 & 8 & 3 & 8 & $7k$ & 4 & $4k$ & $F(72, 3^{96k}4^{7k})$ & $100k,\ k=1,
\ldots, 6$ \\
\hphantom{1}9 & 16 & 3 & 8 & $7k$ & 4 & $2k$ & $F(72, 3^{96k}4^{7k})$ & $100k,\ k=2, \ldots, 6$\\
\hphantom{1}6 & 5 & 3& 12& $11k$ & 6 & $10k$ & $F(72, 3^{150k}6^{11k})$ & $156k,\ k=1, 2, 3$\\
\hphantom{1}8 & 7 & 2 & 10 & $9k$ & 5 & $20k$ & $F(80, 2^{280k}5^{9k})$ & $285k,\ k=1, \ldots, 4$ \\
10 & 18 & 2 & 8 & $7k$ & 4 & $6k$ & $F(80, 2^{216k}4^{7k})$ & $220k,\ k=1, \ldots, 6$\\
\hphantom{1}6 & 5 & 3 & 14 & $13k$ & 7 & $14k$ & $F(84, 3^{210k}7^{13k})$ & $217k,\ k=1, 2$ \\
\hphantom{1}6 & 5 & 3 & 15 & $7k$ & 5 & $5k$ & $F(90, 3^{75k}2^{7k})$ & $80k,\ k=2,
\ldots, 11$\\
\hphantom{1}6 & 10 & 3 & 16 & 30 & 4 & 8 & $F(96, 3^{240}4^{30})$ & 264 \\
10 & 18 & 2 & 10 & $9k$ & 5 & $10k$ & $F(100, 2^{360k}2^{9k})$ &
$365k,\ k=1, \ldots, 4$ \\
\hphantom{1}9 & 4 & 3 & 12 & $11k$ & 2 & $4k$ & $F(108, 3^{48k}2^{11k})$ & $58k,\ k=1, \ldots, 12$ \\
10 & 18 & 2 & 12 & $11k$ & 3 & $6k$ & $F(120, 2^{216k}3^{11k})$ &
$225k,\ k=1, \ldots, 5$ \\
\hphantom{1}8 & 7 & 4 & 14 & $26$ & 7 & $14$ & $F(112, 4^{392}4^{26})$ & $406$\\
\hphantom{1}8 & 14 & 2 & 16 & $5k$ & 4 & $4k$ & $F(128, 2^{112k}4^{5k})$ & $116k,\ k=2, \ldots, 16$ \\
\hphantom{1}9 & 4 & 3 & 15& $7k$ & 5 & $10k$ & $F(135, 3^{120k}5^{7k})$ & $125k,\ k=1, \ldots, 5$\\
\hline
\end{tabular*}
\tabnotetext[\dag]{tab2}{$F(n_1n_2, q^{rm_1q_1}_1q^{m_2}).$}
\legend{$\omega$ is the constant natural weighted
coincidence number of the final resulting SSD.}
\end{table}

\begin{table}
\tabcolsep=0pt
\caption{Some selected $E(f_{\mathrm{NOD}})$-optimal SSDs constructed by
Theorem \protect\ref{sthm4}}
\label{tt3}
\begin{tabular*}{\textwidth}{@{\extracolsep{\fill}}lccccccccccccc@{}}
\hline
$\bolds{n_1}$&$\bolds{m_1}$&$\bolds{q_1}$&$\bolds{n_2}$&$\bolds{m_2}$&$\bolds{q_2}$&$\bolds{m_3}$&$\bolds{q_3}$&$\bolds{m_4}$&$\bolds{q_4}$&
$\bolds{r_3}$ & $\bolds{r_4}$ & \textbf{Final resulting SSD\tabnoteref[\dag]{tab3}} & $\bolds\lambda$\\\hline
\hphantom{1}4 & 3 & 2 & 6 & 10 & 2 & 3 & 2 & 5 & 3 & 6 & 4 & $F(24,
2^{36}3^{60}4^{30})$ & 42 \\
\hphantom{1}4 & 3 & 2 & 6 & 10 & 3 & 3 & 2 & 5 & 3 & 8 & 2 & $F(24,
2^{48}3^{30}6^{30})$ & 36 \\
\hphantom{1}4 & 3 & 2 & 6 & 15 & 3 & 3 & 2 & 5 & 3 & 12 & 3 & $F(24,
2^{72}3^{45}6^{45})$ & 54 \\
\hphantom{1}4 & 3 & 2 & 8 &$7k$ & 2 & 3 & 2 & 7 & 4 & $4k$ & $2k$ & $F(32,
2^{24k}4^{77k})$ & $29k,\ k=1, \ldots, 5$ \\
\hphantom{1}4 & 3 & 2 & 8 & 21 & 2 & 3 & 2 & 14 & 4 & 12 & 3 & $F(32,
2^{72}4^{231})$ & 87 \\
\hphantom{1}4 & 3 & 2 & 8 & 28 & 2 & 3 & 2 & 7 & 4 & 16 & 8 & $F(32,
2^{96}4^{308})$ & 116 \\
\hphantom{1}4 & 3 & 2 & 8 & 21 & 4 & 3 & 2 & 7 & 4 & 18 & 2 & $F(32,
2^{108}4^{56}8^{63})$ & 71 \\
\hphantom{1}6 & 10& 2 & 6 & 10 & 2 & 10& 2 & 5 & 3 & 12 & 12 & $F(36,
2^{240}3^{180}4^{100})$ & 196 \\
\hphantom{1}4 & 3 & 2 & 9 & $8k$ & 3 & 3 & 2 & 4 & 3 & $6k$ & $4k$ & $F(36,
2^{36k}3^{48k}6^{24k})$ & $36k,\ k=1, \ldots, 6$ \\
\hphantom{1}6 & 10& 2 & 6 & $5k$ & 3 & 10& 2 & 5 & 3 & $8k$ & $3k$ & $F(36,
2^{160k}3^{45k}6^{50k})$ & $99k,\ k=1, 2, 3$ \\
\hphantom{1}6 & 10& 2 & 6 & 10 & 3 & 10& 3 & 10 & 2 & 8 & 6 & $F(36,
2^{120}3^{240}6^{100})$ & 148 \\
\hphantom{1}6 & 5 & 3 & 6 & 10 & 2 & 5 & 3 & 10 & 2 & 3 & 8 & $F(36,
2^{160}3^{45}6^{50})$ & 99 \\
\hphantom{1}6 & 5 & 3 & 6 &$5k$ & 3 & 10& 2 & 5 & 3 & $2k$ & $2k$ & $F(36,
2^{40k}3^{30k}9^{25k})$ & $31k,\ k=2, 3$ \\
\hphantom{1}6 & 5 & 3 & 6 & 15 & 3 & 10& 3 & 10 & 2 & 3 & 6 & $F(36,
2^{120}3^{90}9^{75})$ & 93 \\
\hphantom{1}6 & 10& 2 & 8 & 7 & 2 & 10& 3 & 7 & 2 & 4 & 18 & $F(48,
2^{252}3^{120}4^{70})$ & 178 \\
\hphantom{1}6 & 5 & 3 & 8 & 14 & 2 & 5 & 3 & 14 & 2 & 4 & 12 & $F(48,
2^{336}3^{60}6^{70})$ & 194 \\
\hphantom{1}6 & 5 & 3 & 8 & 21 & 2 & 10& 3 & 14 & 2 & 3 & 18 & $F(48,
2^{504}3^{90}6^{105})$ & 291 \\
\hphantom{1}6 & 5 & 3 & 8 & 28 & 2 & 5 & 3 & 21 & 2 & 8 & 16 & $F(48,
2^{672}3^{120}6^{140})$ & 388 \\
\hphantom{1}6 & 5 & 3 & 8 & $7k$ & 4 & 5 & 3 & 7 & 2 & $3k$ & $4k$ & $F(48,
2^{56k}3^{45k}12^{35k})$ & $44k,\ k=1, \ldots, 6$ \\
\hphantom{1}6 & 10& 2 & 8 & $7k$ & 2 & 10& 2 & 14 & 4 & $8k$ & $3k$ & $F(48,
2^{160k}4^{238k})$ & $134k,\ k=1, \ldots, 5$ \\
\hphantom{1}6 & 5 & 3 & 8 & $7k$ & 2 & 10& 2 & 7 & 4 & $2k$ & $4k$ & $F(48,
2^{40k}4^{112k}6^{35k})$ & $ 51k,\ k=2, \ldots, 5$ \\
\hphantom{1}6 & $5k$ & 3 & 8 & 21 & 2 & 10& 2 & $7k$ & 4 & $6k$ & 12 & $F(48,
2^{120k}4^{336k}6^{105k})$ & $153k,\ k=1, 2, 3$ \\
\hphantom{1}6 & 10 & 2 & 8 & $7k$ & 2 & 5 & 3 & 7 & 4 & $8k$ & $6k$ & $F(48,
3^{120k}4^{238k})$ & $94k,\ k=1, \ldots, 5$ \\
\hphantom{1}6 & 5 & 3 & 8 & 14 & 2 & 5 & 3 & 7 & 4 & 4 & 8 & $F(48,
3^{60}4^{224}6^{70})$ & 82 \\
\hphantom{1}6 & $5k$ & 3 & 8 & 21 & 2 & $5k$ & 3 & $7k$ & 4& 6 & 12 & $F(48,
3^{90k}4^{336k}6^{105k})$ & $123k,\ k=1, 2, 3$\\
\hphantom{1}6 & 10& 3 & 9 & 8 & 3 & 10& 2 & 8 & 3 & 6 & 8 & $F(54,
2^{120}3^{192}9^{80})$ & 128 \\
10 & $18$ & 5 & 6 & 10 & 2 & 18 & 2 & 5 & 3 & $6$ & $32$ & $F(60,
2^{216}3^{480}10^{180})$ & $276$\\
\hphantom{1}6 & 10& 3 & 10 & 9 & 5 & 5 & 3 & 18 & 2 & 8 & 4 & $F(60,
2^{144}3^{120}15^{90})$ & 114 \\
\hphantom{1}6 & 10& 3 & 10 & 9 & 5 & 5 & 3 & 9 & 5 & 8 & 2 & $F(60,
3^{120}5^{90}15^{90})$ & 60 \\
\hphantom{1}8 & 7 & 2 & 8 & $7k$ & 2 & 7 & 2 & 7 & 4 & $12k$ & $4k$ & $F(64,
2^{168k}4^{161k})$ &$121k,\ k=1, \ldots, 5$\\
\hphantom{1}8 & $7k$ & 2 & 8 & 14 & 2 & 14 & 2 & 7 & 4 & $12k$ & $8k$ & $F(64,
2^{336k}4^{322k})$ &$242k,\ k=1, \ldots, 5$\\
\hphantom{1}8 & 7 & 2 & 8 &$7k$ & 4 & 21 &4 & 7 & 2 &$2k$ & $4k$ & $F(64,
2^{56k}4^{168k}8^{49k})$ &$73k,\ k=1, \ldots, 6$ \\
\hphantom{1}8 & $7k$ & 2 & 8 &14 & 4 & 21 &4 & 14 & 2 &$4k$ & $4k$ & $F(64,
2^{112k}4^{336k}8^{98k})$ & $146k,\ k=1, \ldots, 5$ \\
\hphantom{1}8 & $7k$ & 2 & 8 &21 & 4 & 21 &2 & 7 & 4 &$18k$ & $4k$ & $F(64,
2^{756k}4^{112k}8^{147k})$ &$415k,\ k=1, \ldots, 5$ \\
\hphantom{1}8 & $7k$ & 2 & 9 & 16 & 3 & $7k$ &4 & 16 & 3 &12 & $4k$ & $F(72,
3^{192k}4^{336k}6^{112k})$ &$160k,\ k=1, \ldots, 5$ \\
10 & $9k$ & 5 & 8 & 14 & 2 & $9k$ & 5 & $7k$ & 4 & 2& 16 & $F(80,
4^{448k}5^{90k}10^{126k})$ & $136k,\ k=1, \ldots, 4$\\
10 & $9k$ & 5 & 8 & 21 & 2 & $9k$ & 5 & $14k$ & 4 & 3 & $12$ & $F(80,
4^{672k}5^{135k}10^{189k})$& $204k,\ k=1, 2, 3$\\
10 & $9k$ & 5 & 8 & 14 & 4& $9k$ & 5 & $7k$ & 2 & 3 & 16 & $F(80,
2^{224k}5^{135k}20^{126k})$ & $141k,\ k=1, \ldots, 4$\\
\hline
\end{tabular*}
\tabnotetext[\dag]{tab3}{$F(n_1n_2,
(q_1q_2)^{m_1m_2}q^{m_3r_3q_3}_3q^{m_4r_4q_4}_4)$; $n_1=n_3$,
$n_2=n_4$.}
\legend{$\lambda$ is the constant coincidence number of
the final resulting SSD.}
\end{table}

\begin{table}
\tabcolsep=0pt
\caption{Some selected $\chi^2$-optimal SSDs constructed by Theorem \protect\ref{sthm4}}
\label{tt2}
\begin{tabular*}{\textwidth}{@{\extracolsep{\fill}}lccd{2.0}cccccccccc@{}}
\hline
$\bolds{n_1}$&$\bolds{m_1}$&$\bolds{q_1}$&\multicolumn{1}{c}{$\bolds{n_2}$}&$\bolds{m_2}$&$\bolds{q_2}$&$\bolds{m_3}$
&$\bolds{q_3}$&$\bolds{m_4}$&$\bolds{q_4}$&
$\bolds{r_3}$ & $\bolds{r_4}$ & \textbf{Final resulting SSD\tabnoteref[\dag]{tab4}} & $\bolds\omega$\\ \hline
\hphantom{1}4 & 3 & 2 & 6 & 5 & 3 & 3 & 2 & 5 & 3& 12 & 2 & $F(24,
2^{72}3^{30}6^{15})$ & 108 \\
\hphantom{1}4 & 3 & 2 & 8 & $7k$ & 2 & 3 & 2 & 7 & 4& $8k$ & $2k$ & $F(32,
2^{48k}4^{77k})$ & $116k,\ k=1, \ldots, 5$ \\
\hphantom{1}9 & $8k$ & 3 & 4 & 3 & 2& $8k$ & 3 & 3 & 2 & 4 & $18k$ & $F(36,
2^{108k}3^{96k}6^{24k})$ & $216k,\ k=1, \ldots, 6$ \\
\hphantom{1}6 & $5k$ & 3 & 6 & 10 & 3 & 10& 2 & 10 & 3 & $18k$ & $6k$ & $F(36,
2^{360k}3^{180k}9^{50k})$ & $558k,\ k=1, 2, 3$\\
\hphantom{1}6 & 10 & 2 & 8 & 7 & 2 & 10& 2 & 14& 4& 16 & 3 & $F(48,
2^{320}4^{238})$ & 536 \\
\hphantom{1}6 & 5 & 3 & 8 & 7 & 2 & 10& 2 & 14& 4& 6 & 3 & $F(48,
2^{120}4^{168}6^{35})$ & 306 \\
\hphantom{1}6 & 5 & 3 & 8 & 7 & 2 & 5 & 3 & 14& 2& 4 & 18 & $F(48,
2^{504}3^{60}6^{35})$ & 582 \\
\hphantom{1}6 & 5 & 3 & 8 & $7k$ & 2 & 5 & 3 & $14k$ & 4& $4k$ & 3 & $F(48,
3^{60k}4^{168k}6^{35k})$ & $246k,\ k=1, 2, 3$ \\
\hphantom{1}6 & 5 & 3 & 8 & 14& 2 & 10& 2 & 7 & 4& 12 & 12 & $F(48,
2^{240}4^{336}6^{70})$ & 612 \\
\hphantom{1}6 & 5 & 3 & 8 & 14& 2 & 5 & 3 & 7 & 4& 8 & 12 & $F(48,
3^{120}4^{336}6^{70})$ & 492 \\
\hphantom{1}6 & 5 & 3 & 8 & 21& 2 & 10& 2 & 42& 4& 18 & 3 & $F(48,
2^{360}4^{504}6^{105})$ & 918 \\
\hphantom{1}6 & 5 & 3 & 8 & 28& 2 & 10& 3 & 42& 4& 8 & 4 & $F(48,
3^{240}4^{672}6^{140})$ & 984 \\
\hphantom{1}6 & 5 & 3 & 8 & 7 & 4 & 10& 2 & 7 & 4& 18 & 4 & $F(48,
2^{360}4^{112}12^{35})$ & 484 \\
\hphantom{1}6 & 5 & 3 & 8 & 7 & 4 & 5 & 3 & 14& 2& 12 & 12 & $F(48,
2^{336}3^{180}12^{35})$ & 528 \\
\hphantom{1}6 & 5 & 3 & 8 & 7 & 4 & 5 & 3 & 7 & 4& 12 & 4 & $F(48,
3^{180}4^{112}12^{35})$ & 304 \\
\hphantom{1}6 & 5 & 3 & 8 & 14& 4 & 10& 3 & 21& 2& 12 & 16 & $F(48,
2^{672}3^{360}12^{70})$ & 1056 \\
\hphantom{1}6 & 5 & 3 & 8 & $7k$& 4 & 10& 3 & 7 & 4& $6k$ & $4k$ & $F(48,
3^{180k}4^{112k}12^{35k})$ & $304k,\ k=1, \ldots, 6$ \\
\hphantom{1}6 & 5 & 3 & 8 & 21& 4 & 10& 3 & 7 & 4& 18 & 12 & $F(48,
3^{540}4^{336}12^{105})$ & 912 \\
\hphantom{1}6 & 10 & 2 & 9 & 4 & 3 & 10& 2 & 4 & 3& 18 & 12 & $F(54,
2^{360}3^{144}6^{40})$ & 528 \\
10 & $9k$ & 5 & 6 & 5& 3 & $9k$ & 5 & $5k$ & 3 & 3 & 20 & $F(60,
3^{300k}5^{135k}15^{45k})$ & $450k,\ k=1, 2, 3$\\
\hphantom{1}8 & 7 & 2 & 8 & 7 & 2 & 14& 2 & 7 & 4& 12 & 4 & $F(64,
2^{336}4^{161})$ & 484 \\
\hphantom{1}8 & 7 & 2 & 8 & 14& 2 & 21& 2 & 7 & 4& 16 & 8 & $F(64,
2^{672}4^{322})$ & 968 \\
\hphantom{1}8 & 7 & 2 & 8 & 21& 2 & 28& 2 & 7 & 4& 18 & 12 & $F(64,
2^{1008}4^{483})$ & 1452 \\
\hphantom{1}8 & 7 & 2 & 8 & $7k$& 4 & $7k$ &4 & $7k$ & 2& 12 & 16 & $F(64,
2^{224k}4^{336k}8^{49k})$ & $584k,\ k=1, \ldots, 5$ \\
\hphantom{1}8 & 7 & 2 & 8 & 21& 4 & 21& 4 & 28& 2& 12 & 12 & $F(64,
2^{672}4^{1008}8^{147})$ & 1752 \\
\hphantom{1}8 & 7 & 2 & 9 & 8 & 3 & 21& 2 & 8 & 3& 18 & 8 & $F(72,
2^{756}3^{192}6^{56})$ & 984 \\
\hphantom{1}8 & 7 & 2 & 9 &$8k$& 3 & $21k$& 4& $16$ & 3& 3 & $4k$ & $F(72,
3^{192k}4^{252k}6^{56k})$ & $480k,\ k=1, 2$ \\
\hphantom{1}8 & 7 & 4 & 9 & 4 & 3 & 7 & 2 & 8 & 3& 18 & 12 & $F(72,
2^{252}3^{288}12^{28})$& 552 \\
\hphantom{1}8 & 7 & 4 & 9 & $4k$& 3 & $7k$& 4& $8k$ & 3& 3 & 12 & $F(72,
3^{288k}4^{84k}12^{28k})$& $384k,\ k=1, \ldots, 6$ \\
\hphantom{1}9 & $8k$ & 3 & 8 & 7 & 2& $8k$ & 3 & 21 & 2 & 8 & $18k$ & $F(72,
2^{756k}3^{192k}6^{56k})$ & $984k,\ k=1, \ldots, 6$ \\
\hphantom{1}9 & $8k$ & 3 & 8 & 7 & 2& $8k$ & 3 & 21 & 4 & 8 & $3k$ & $F(72,
3^{192k}4^{252k}6^{56k})$ & $480k,\ k=1, \ldots, 6$ \\
\hphantom{1}9 & $8k$ & 3 & 8 & 7 & 4& 16 & 3 & 21 & 2 & $12k$ & $12k$ & $F(72,
2^{504k}3^{576k}12^{56k})$ &$1104k,\ k=1, \ldots, 4$ \\
\hphantom{1}9 & $8k$ & 3 & 8 & 7 & 4& 16 & 3 & 14 & 4 & $12k$ & $3k$ & $F(72,
3^{576k}4^{168k}12^{56k})$ & $768k,\ k=1, \ldots, 4$ \\
10 & $9$ & 5 & 8 & $7k$ & 2 & 18 & 2 & $7k$ & 4 & $10k$ & 20 & $F(80,
2^{360k}4^{560k}10^{63k})$& $950k,\ k=1, \ldots, 5$\\
\hphantom{1}8 & 14 & 2 & 10 & $9k$ & 5 & 28 &2 & $9k$ & 5& $12k$ & 2 & $F(80,
2^{672k}5^{90k}10^{126k})$ &$360k,\ k=1, \ldots, 4$\\
10 & $9k$ & 5 & 8 & 7 & 4 & $9k$ & 5 & 28 & 2 & 6 & $20k$ & $F(80,
2^{1120k}5^{270k}20^{63k})$ & $1410k,\ k=1, \ldots, 5$\\
\hphantom{1}8& $7k$ & 2& 10 & 9 & 5 & $7k$ & 4 & 9 & 5 & 20 & $2k$ & $F(80,
4^{560k}5^{90k}10^{63k})$ & $680k,\ k=1, \ldots, 5$\\
\hline
\end{tabular*}
\tabnotetext[\dag]{tab4}{$F(n_1n_2,
(q_1q_2)^{m_1m_2}q^{m_3r_3q_3}_3q^{m_4r_4q_4}_4)$;
$n_1=n_3$, $n_2=n_4$.}
\legend{$\omega$ is the constant natural weighted
coincidence number of the final resulting SSD.}
\vspace{3pt}
\end{table}

\end{appendix}

\newpage

\section*{Acknowledgments}
The authors thank the Editor, the Associate Editor
and two referees for their valuable and constructive comments which
have led to
a significant improvement in the presentation of this paper.


%

\printaddresses

\end{document}